\documentclass[ijoc]{informs4}
\usepackage{graphicx} % Required for inserting images
\usepackage{amsmath,amssymb}
\newcommand{\dotsint}{\mathinner{{\ldotp}{\ldotp}}}

\usepackage{subcaption}
\usepackage{xcolor}
\usepackage{multirow}
\usepackage{placeins}
\usepackage{booktabs}

\usepackage{graphicx}
\usepackage[export]{adjustbox}

\usepackage[ruled,vlined, linesnumbered]{algorithm2e}
\SetCommentSty{textrm}
\SetAlFnt{\small}
\SetKwInput{Input}{Input}
\SetKwInput{Output}{Input}
\SetArgSty{textnormal}

\usepackage{newrevisor}
\newrevisor{manuel}{orange}
\newrevisor{isaac}{purple}

\newcommand{\abs}[1]{\ensuremath{\lvert#1\rvert}}
\newcommand{\norm}[1]{\ensuremath{\lVert#1\rVert_2}}
\newcommand{\assign}{\ensuremath{\leftarrow}}
\newcommand{\tmax}{\ensuremath{t_\text{max}}}
\newcommand{\taumax}{\ensuremath{\tau_\text{max}}}
\newcommand{\tauf}{\ensuremath{\tau_\text{f}}}
\newcommand{\Earth}{\textit{Earth}\xspace}
\newcommand{\Done}{\textit{done}\xspace}

\newcommand{\earlieststart}{\ensuremath{\mathit{est}}}
\newcommand{\earliestarrival}{\ensuremath{\mathit{eat}}}
\newcommand{\Btilde}{\ensuremath{\tilde{\mathcal{B}}}}
\newcommand{\zopt}{\ensuremath{z_\text{opt}}}
\newcommand{\Bprime}{\ensuremath{\mathcal{B}^\prime}}
\newcommand{\nodelabel}{\ensuremath{l}}
\newcommand{\arclabel}{\ensuremath{l}}

\usepackage{hyperref}
\hypersetup{colorlinks,allcolors=blue}
\usepackage{natbib}
 \bibpunct[, ]{(}{)}{,}{a}{}{,}%
 \def\bibfont{\small}%
\JOURNAL{Informs Journal on Computing}
\OneAndAHalfSpacedXI
\TheoremsNumberedThrough  
\ECRepeatTheorems
\EquationsNumberedThrough  

\begin{document}
\RUNTITLE{Exact Solutions to the Space-Time Dependent TSP}

\TITLE{An Exact Framework for Solving the Space-Time Dependent TSP}

\ARTICLEAUTHORS{%
\AUTHOR{Isaac Rudich}
\AFF{Polytechnique Montréal, Montreal, Canada, \EMAIL{isaac.rudich@polymtl.ca}}

\AUTHOR{Manuel López-Ibáñez}
\AFF{University of Manchester, Manchester, England, \EMAIL{manuel.lopez-ibanez@manchester.ac.uk}}

\AUTHOR{Michael Römer}
\AFF{Universität Bielefeld, Bielefeld, Germany, \EMAIL{michael.roemer@uni-bielefeld.de}}

\AUTHOR{Quentin Cappart, Louis-Martin Rousseau}
\AFF{Polytechnique Montréal, Montreal, Canada, \{quentin.cappart@polymtl.ca, louis-martin.rousseau@polymtl.ca\}}
}

\ABSTRACT{
    Many real-world scenarios involve solving bi-level optimization problems in which there is an outer discrete optimization problem, and an inner problem involving expensive or black-box computation. This arises in space-time dependent variants of the Traveling Salesman Problem, such as when planning space missions that visit multiple astronomical objects. Planning these missions presents significant challenges due to the constant relative motion of the objects involved. There is an outer combinatorial problem of finding the optimal order to visit the objects and an inner optimization problem that requires finding the optimal departure time and trajectory to travel between each pair of objects. The constant motion of the objects complicates the inner problem, making it computationally expensive.
    This paper introduces a novel framework utilizing decision diagrams (DDs) and a DD-based branch-and-bound technique, Peel-and-Bound, to achieve exact solutions for such bi-level optimization problems, assuming sufficient inner problem optimizer quality. The framework leverages problem-specific knowledge to expedite search processes and minimize the number of expensive evaluations required.
    As a case study, we apply this framework to the Asteroid Routing Problem (ARP), a benchmark problem in global trajectory optimization. Experimental results demonstrate the framework's scalability and ability to generate robust heuristic solutions for ARP instances. Many of these solutions are exact, contingent on the assumed quality of the inner problem's optimizer. 
}

\MSCCLASS{90-08}

\maketitle

% \listofrevisions

\section{Introduction}
In various real-world scenarios, we must determine the optimal sequence for visiting several locations or scheduling jobs, where the travel cost between two locations, or the setup time between jobs, is defined by an expensive black-box function without a closed form. These scenarios include a class of bi-level optimization problems. The \emph{outer} problem involves finding an optimal permutation, while evaluating the cost of a given, possibly partial, permutation requires solving a different \emph{inner} optimization problem that depends on the permutation. When the inner problem is computationally expensive and/or black-box, finding exact solutions to the bi-level optimization problem is difficult in practice.

Instances of this class of problems include time-dependent routing problems~\citep{GenGhiGue2015tdtsp}, especially variants where the travel costs must be computed on-demand \citep{EhmCamTho2016tdvrp} or are given by solving an inner optimization problem that, for example, determines vehicle speed \citep{AndFagHob2015maritime}. This problem class also arises when planning missions in Space \citep{ShiCebLoz2018space}, such as finding a tour of the Jupiter Galilean moons \citep{IzzSimMar2013tour}, global optimization for multiple gravity assist trajectories \citep{AbdGad2012dynamic,CerVas2010mga,IzzBecMyaNas2007search,VasPas2006preliminary}, active space debris removal \citep{IzzGetHenSim2015evolving} and spacecraft formation control \citep{TaeLeoClam2007spacecraft}. These problems are often solved through bi-level heuristics that combine numerical optimization methods with tree search methods such as Beam-Search \citep{PetBonGre2014gtoc5}, Lazy Race Tree Search \citep{IzzSimMar2013tour}, Monte-Carlo Tree Search \citep{HenIzz2015interplanetary} and Beam-ACO \citep{SimIzzHaas2017multi}. Other approaches treat the inner problem as a black-box function whose evaluation is expensive \citep{ZaeStoFriFisNauBar2014,IruLop2021gecco,SanBai2022permutations,ChiDerVer2023fourier}. To the best of our knowledge, no exact method has ever been proposed to find optimal solutions for such problems. In the context of Space exploration, solutions to global trajectory optimization problems represent huge economic costs, hundreds of millions of dollars, and long mission horizons, from years to decades. Thus, optimality is highly desirable.

We propose a framework for finding exact solutions by representing the outer problem using decision diagrams (DDs) \citep{BerCirHoeHoo2016dd4o,CirHoe2013mdd} and then solving that model using a DD-based branch-and-bound technique called Peel-and-Bound \citep{RudCapRou2023improved_pnb} with DD-based search techniques \citep{Gillard2022phd}. DDs are a tool for using graphs to compactly represent the solution space of discrete optimization problems \citep{BerCirHoeHoo2016dd4o}. Our proposed approach treats the inner problem as an expensive black-box function evaluated on-demand. Still, it uses problem-specific knowledge about this black-box function to speed-up the search, and reduce the number of black-box evaluations needed to solve the problem. 

As a case study, we consider the Asteroid Routing Problem (ARP) \citep{LopChiGil2022evo}, which is a benchmark problem inspired by the $11^\text{th}$ Global Trajectory Optimization Competition (\url{https://gtoc11.nudt.edu.cn}). The ARP is a space-time dependent variant of the TSP. In the ARP, a spacecraft launched from Earth is tasked with visiting a specific set of asteroids. The goal is to find the permutation of asteroids that minimizes fuel consumption and total travel time, aggregated into a single objective function. Traveling from one asteroid to another requires identifying optimal departure and travel times, which constitutes the inner problem. The ARP treats this inner problem as a black-box function that a pre-defined deterministic optimizer solves. Thus, in the ARP, evaluating a given permutation of the asteroids always yields the same objective function value. These characteristics make the ARP a simplified benchmark for the class of problems described above and allow researchers to focus on the \textit{outer} problem of optimizing the permutation of asteroids. As in the problems mentioned above, a brute force search is the only known method for identifying an optimal solution to the ARP, which is only computationally tractable for a very small number of asteroids. All approaches proposed so far are heuristic and time-consuming \citep{LopChiGil2022evo,ChiDerVer2023fourier}.

This paper proposes the first exact approach for solving the ARP. Our framework proves the optimality of the outer problem under the assumption that the inner problem is optimized with sufficient quality. We address the quality of the inner optimizer in our experimental results (Section \ref{sec:exps}). We also discuss how this technique can be made scalable and used to quickly generate strong heuristic solutions to the ARP. We provide exact solutions and new best-known solutions for several ARP instances.

This paper is structured as follows. Section \ref{sec:background} provides the necessary background and notation for the ARP and DD-based methods. Section \ref{sec:initialconstruction} details the construction of the initial relaxed DD. Section \ref{sec:heuristicsearch} outlines how to use embedded restricted DDs within the relaxed DD to search for solutions. Section \ref{sec:pnb} explains how Peel-and-Bound can be used on the relaxed DD to find exact solutions. Section \ref{sec:implementation} details key implementation decisions and insights. Section \ref{sec:exps} provides our experimental results. Finally, Sections \ref{sec:future} and \ref{sec:conclusion} summarize our conclusions and suggests directions for further research.

\section{Background}
\label{sec:background}
In this section we detail the ARP, the relevant DD techniques, and how they apply to the ARP.
% We fist define the general black-box time-dependent routing problem. Next, we define the ARP and explain how the ARP is an example of such problems.

% \subsection{Black-Box Time-Dependent Routing Problems}

% We consider a class of routing problems where we must find an optimal permutation $\pi$ of the set of $n$ locations $A=\{a_1, \dotsc, a_n\}$ that minimizes a function $f\colon S_n \to \mathcal{R}$, where $S_n$ is the set of all permutations of $A$. The value of $f(\pi)$ is given by $\sum_[i=2]^n \mathcal(Q)(\pi_{i-1}, \pi_{i}, t_i)$, where $\pi_i \in A$ and $t_i$ is a parameter that depends on when $\pi_i$ is visited and, thus, it also depends on the partial permutation $\{pi_1, \dotsc, \pi_i\}$. If a closed-form of $\mathcal{B}$ is given, then this formulation corresponds to the well-known time-dependent TSP problem. If evaluating $\mathcal{B}$ requires solving an optimization problem 

\subsection{The Asteroid Routing Problem}
\label{sec:ARP}
In the ARP \citep{LopChiGil2022evo}, we are given a set of $n$ asteroids \mbox{$A = \{a_1,\dotsc,a_n\}$} orbiting around the Sun that a spacecraft launched from Earth ($a_0$) must visit. The spacecraft does not need to return to Earth. For simplicity, the ARP treats the spacecraft transfer from Earth to the first-visited asteroid as any other transfer between asteroids, i.e., Earth has no escape velocity.

\begin{figure}[h!t]
    \FIGURE{
        \label{fig:arp_drawing}
        \centering
        \begin{tabular}{c@{\hspace{3em}}c}
            \begin{subfigure}[t]{0.45\linewidth}
                \centering
                \includegraphics[scale=.8]{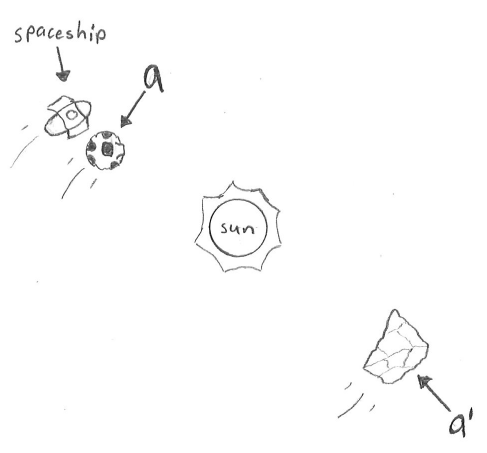}
                \captionsetup{font=normal}
                \caption{At time $\eta$, the spaceship in in $a$'s orbit.}
            \end{subfigure}
            &
            \begin{subfigure}[t]{0.45\linewidth}
                \centering
                \includegraphics[scale=.8]{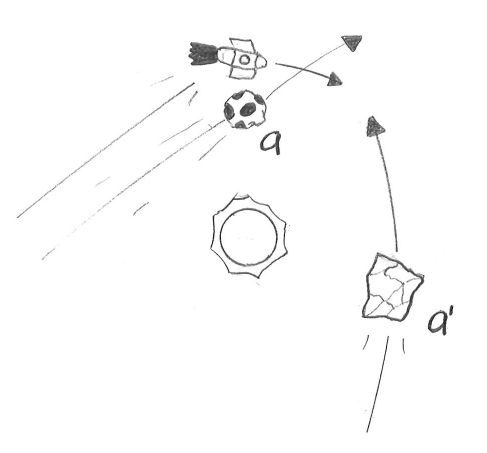}
                \vspace{.92em}
                \captionsetup{font=normal}
                \caption{After $\tau$ days, at time $\eta + \tau$, the spaceship initiates a transfer to $a^\prime$.}
            \end{subfigure} \\
            \begin{subfigure}[t]{0.45\linewidth}
                \centering
                \vspace{1em}
                \includegraphics[scale=.8]{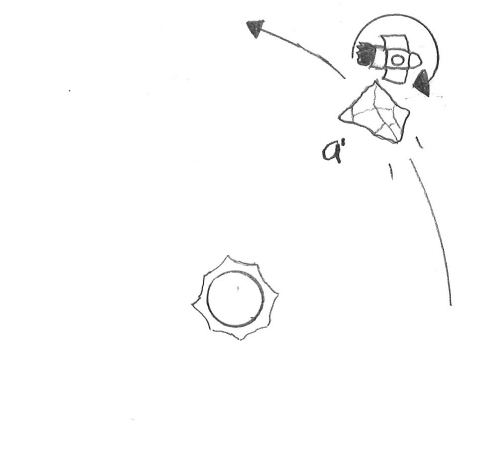}
                \vspace{-2em}
                \captionsetup{font=normal}
                \caption{After $t$ days, at time $\eta + \tau + t$, the spaceship intersects $a^\prime$ and performs a maneuver to align itself with $a^\prime$ 's orbit.}
            \end{subfigure}
            & 
            \begin{subfigure}[t]{0.45\linewidth}
                \centering
                \vspace{1em}
                \includegraphics[scale=.8]{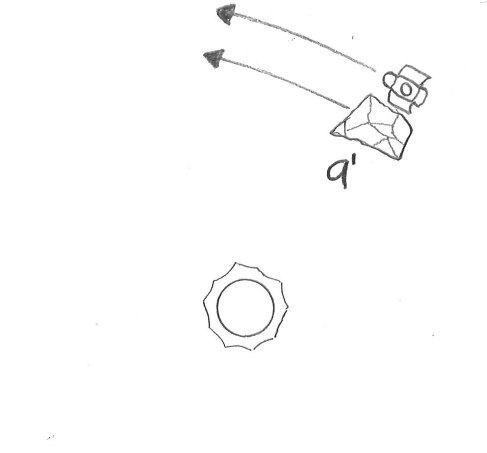}
                \vspace{-2em}
                \captionsetup{font=normal}
                \caption{Still at time $\eta + \tau + t$, the spaceship is now in the same orbit as $a^\prime$.}
            \end{subfigure}
        \end{tabular}
    }
    {Depiction of a transfer from $a$ to $a^\prime$}
    {}
\end{figure}

A transfer of the spacecraft between two asteroids $a$ and $a'$ that starts no earlier than time (epoch) $\eta$ is computed with the following black-box function:
\begin{equation}\label{eq:transfer}
  % \mathcal{B}(a, a', \eta) \to (\tau, t, z)
  \mathcal{B}(a, a', \eta)\colon A \times A \times \mathbb{R}^{+} \to [0, \taumax] \times [1, \tmax] \times \mathbb{R}^{+}
\end{equation}
which returns three values $(\tau, t,z)$, where $\tau$ is how long the spacecraft waits at $a$ before departing, $t$ is how long the spacecraft travels until it arrives at $a'$, such that arrival time is $\eta+\tau+t$, and $z$ is the cost of this transfer. A transfer is depicted in Figure \ref{fig:arp_drawing}.

% \begin{figure}[h!t]
%     \FIGURE{
%         \label{fig:arp_drawing}
%         \centering
%         \includegraphics[scale=.2,frame]{Figures/ARP_Fig.jpg}
%     }
%     {Depiction of a transfer from $a$ to $a^\prime$}
%     {}
% \end{figure}

As shown above, the wait time $\tau$ and travel time $t$ are bounded by $[0,\taumax]$ and $[1,\tmax]$, respectively. The wait time $\tau$ has a lower bound of $0$ because it can depart immediately without waiting. The travel time $t$ has a lower bound of 1 because a transfer will always require a nonzero time (at least 1 day in the ARP). The upper bounds $\taumax$ and $\tmax$ are set to 730 days in the original definition of the ARP \citep{LopChiGil2022evo}.

\begin{figure}[t]
    \FIGURE{
        \label{fig:arp_solution}
        \centering
        \includegraphics[width=\textwidth]{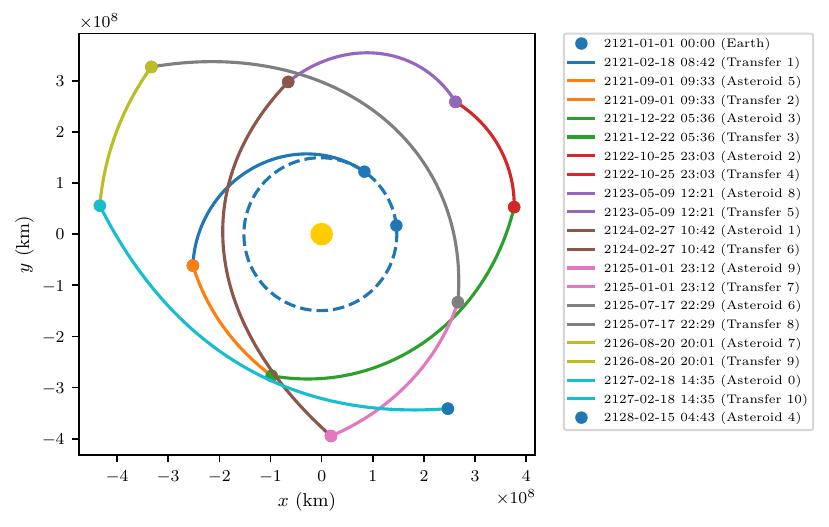}
    }
    {Visualization of best known solution to instance with $n=10$ and \emph{seed}$=8$.}
    {}
\end{figure}

Given a solution $\pi = (\pi_0 = a_0, \pi_1, \dotsc, \pi_n)$, where the subsequence $(\pi_1,\dotsc,\pi_n)$ is a permutation of the asteroids in  $A$ that indicates the order in which they  will be visited, the objective function of the ARP is:
\begin{equation}
  \label{eq:objf}
\text{Minimize}\quad  f(\pi) = \sum_{i=1}^n z_i\qquad\text{where}\quad (\tau_i, t_i, z_i) = \mathcal{B}(\pi_{i-1}, \pi_i, \eta_i)
\end{equation}
Besides,  $\eta_i = \eta_{i-1} + \tau_{i-1} + t_{i-1}$ is the arrival time at $\pi_{i-1}$ and $\eta_1$ is a given mission start time, which can be simplified to be zero. A solution is visualized in Figure \ref{fig:arp_solution}.

The black-box function $\mathcal{B}$ depends not only on the two asteroids being visited but also on the arrival time $\eta_i$. Thus the ARP is an example of the class of time-dependent problems discussed in the introduction.

\subsection{Trajectory Optimization}

In the ARP, and in many trajectory optimization problems, the spacecraft only uses \textit{impulsive} maneuvers. An impulsive maneuver changes the spacecraft's velocity in a negligible amount of time, simplifying the calculation of the maneuver's effect on the spacecraft's orbit. As a result, the maneuver's effect can be modeled as an instant change in the spacecraft's velocity vector $\Delta\vec{v}$, without accounting for the dynamics of the propulsion system or the influence of external forces during the thrust application.

Lambert's Problem \citep{Izzo2015lambert} calculates an orbit that connects two points in space-time departing at time $\eta+\tau$ and arriving at time $\eta+\tau+t$. Using the solution to Lambert's Problem, a transfer between the orbits of two asteroids $a$ and $a'$ can be achieved with two impulsive maneuvers. Initially, at time (epoch) $\eta$, the spacecraft follows the same orbit around the Sun as asteroid $a$. \emph{Waiting} in this orbit does not consume any fuel but does change its relative distance to other asteroids. The first impulse $\Delta\vec{v}_1$ happens $\tau$ days after the current time $\eta$ (epoch), and moves the spacecraft from its current orbit to an orbit that will intercept the target asteroid $a'$ after traveling for $t$ days. When the spacecraft intercepts $a'$, the second impulse $\Delta\vec{v}_2$ modifies its orbit to match the orbit of $a'$, so that the spacecraft follows the same orbit as $a'$ without consuming any additional fuel. Usually, the total magnitude of velocity change (the Euclidean L2 norm) is used as a surrogate for fuel consumption/energy costs:
\begin{equation}
  \Delta V = \norm{\Delta\vec{v}_1} + \norm{\Delta\vec{v}_2}\qquad\text{where}\quad (\Delta\vec{v}_1,\Delta\vec{v}_2) = \text{Lambert}(a,a',\eta+\tau,t)
\end{equation}

We still need to decide $\tau$, how long should the spacecraft wait at $a$, and $t$, how long the transfer should take until arriving at $a'$. Waiting changes the relative distance to other asteroids, which may reduce the fuel or total time needed to reach the next asteroid. The optimal values of $\tau$ and $t$  depend on which objectives are optimized.

In trajectory optimization, there are many possible objectives. Two typical objectives are the minimization of fuel (energy) consumption and the total mission time ($\tau + t$). These two objectives are often in conflict. The ARP aggregates the two objectives as shown below:
\begin{equation}\label{eq:inner}
  \begin{split}
    z = f_\text{inner}(a, a', \eta, \tau, t) &= \Delta V + \frac{2\,\text{km/s}}{30\,\text{days}} \cdot (\tau + t)\\
    \text{s.t.}&\quad \tau \in [0, \taumax],\; t \in [1,\taumax]\\
      \text{where}&\quad\Delta V = \norm{\Delta\vec{v}_1} + \norm{\Delta\vec{v}_2}\quad\text{and}\quad (\Delta\vec{v}_1,\Delta\vec{v}_2) =
\text{Lambert}(a,a',\eta+\tau,t)
  \end{split}
\end{equation}
where the trade-off constant $\frac{2\,\text{km/s}}{30\,\text{days}}$ was chosen by \citet{LopChiGil2022evo}.

The minimization of $f_\text{inner}$ becomes an \emph{inner} optimization problem whose solution gives the values $\tau$, $t$ and $z$ returned by $\mathcal{B}$ (Eq.~\ref{eq:objf}).

In the ARP, this inner problem is solved using Sequential Least Squares Programming (SLSQP) \citep{Kraft1988slsqp}, a deterministic optimizer. In this way, each permutation of the asteroids corresponds to a unique solution. SLSQP iteratively calculates the position and velocity of the asteroids and then
solves Lambert's Problem to evaluate the above objective. Even if time is discretized into days, this inner optimization would be too slow to allow us to calculate the cost of every optimal trajectory for every pair of asteroids at every point in time. 

A mission with $n$ asteroids may have up to $n(\taumax + \tmax)$ possible values of $\eta$. In a mission where $n = 10$, calculating every possible trajectory from just one asteroid to one other asteroid for every possible departure time when time is discretized into days would take about 10 minutes using the inner optimizer in our implementation. Thus, calculating the full cost matrix for every ordered pair of asteroids ($n(n-1)$ pairs) at every possible departure day would take about 15 hours. Furthermore, departure time is continuous, not discrete, so this cost matrix would not even fully represent the problem, although it would certainly provide information that is useful for solving it.

\subsection{Relaxed Decision Diagrams}

\begin{figure}[!tb]
    \FIGURE{
        \centering
        \begin{subfigure}[b]{0.45\linewidth}
            \centering
            \vspace{10pt}
            \includegraphics[scale=.8]{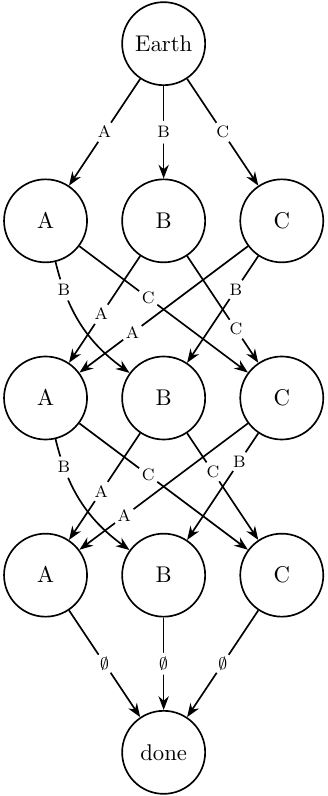}
            \captionsetup{font=normal}
            \caption{Unweighted}
            \label{fig:relaxed_dd_1}
        \end{subfigure}
        \hfill % Optional space between the figures
        \begin{subfigure}[b]{0.45\linewidth}
            \centering
             \vspace{10pt}
            \includegraphics[scale=.8]{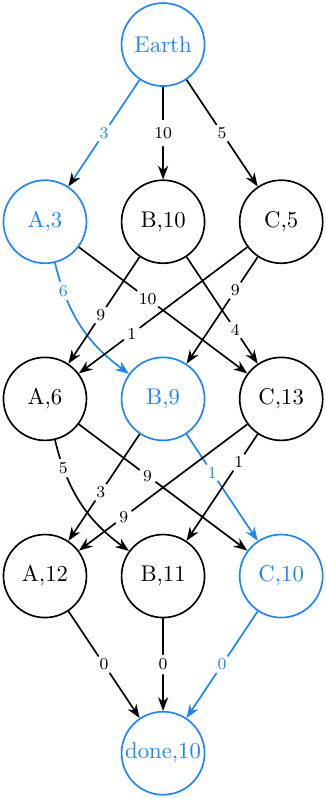}
            \captionsetup{font=normal}
            \caption{Weighted}
            \label{fig:relaxed_dd_2}
        \end{subfigure}
    }
    {Relaxed decision diagrams for an ARP with asteroids  $\{A,B,C\}$.}
    {}
\end{figure}

A decision diagram (DD) is a directed layered graph that encodes solutions to an optimization problem as paths from the root (denoted here as \Earth) to the terminal (denoted here as \Done). A DD is \textit{relaxed} if it encodes every feasible solution, but also encodes infeasible solutions~\citep{CirHoe2013mdd}. Figure~\ref{fig:relaxed_dd_1} displays a valid relaxed DD for an ARP with asteroids $\{A,B,C\}$. Each arc is labeled with the decision about which asteroid is the next destination, and each node is labeled with the current asteroid. Every feasible permutation of asteroids is trivially encoded as a path from the root to the terminal, but so are several infeasible solutions, such as $A \rightarrow B \rightarrow A$.

A relaxed DD can also be \textit{weighted}. A weighted DD is one where every arc is bound by the impact on the objective function of transitioning from the arc's origin to its destination. In weighted relaxed DDs, the length of the shortest path from the root to the terminal is a relaxed (dual) bound on the optimization problem being represented. Figure~\ref{fig:relaxed_dd_2} displays a weighted version of Figure~\ref{fig:relaxed_dd_1}. Each arc still assigns the same decision, but the arc weights now display the lower bounds on the cost of transferring between those two asteroids at that point in the sequence. Each node is now labeled with the shortest path length to that node. Thus, the length of the shortest path ($A \rightarrow B \rightarrow C$, which is highlighted) is a valid lower bound on this ARP. In other words, no matter what permutation of asteroids is used, and how well the trajectories between them are optimized, there is no route for the spacecraft that admits a cost of less than $10$ for this example.

In the well-known Travelling Salesman Problem (TSP), the weights on the arcs would be the transition costs of going from one city to the next. However, in the TSP, the arc weights are not dependent on all preceding cities, and the transition costs are easy to calculate. In this paper, we are handling a problem where calculating arc weights requires calling a black-box function $\mathcal{B}$ with non-negligible computation time. Furthermore, the exact arc weight for arc $a$ depends on the path taken to reach $a$. In principle, this dependency requires the evaluation of all possible paths from $r$ to $a$ to find the optimal one. However, if we can find, via relaxation of the $\mathcal{B}$ function, a globally valid bound on the arc weights, that bound can be used to prune sub-optimal paths.

In the case of the ARP, this relaxation is given by a no-wait variant of $\mathcal{B}$ where the cost of waiting is removed from the objective function; this is defined in Section~\ref{sec:qrelaxed}. This no-wait variant can be used to find a globally valid bound on the cost of any transition between two asteroids. Thus for the ARP, the arcs will be weighted with valid lower bounds on their transition costs instead of exact values. This idea is explored and formalized in Section \ref{sec:initialconstruction}.
   
\subsection{Refining a Relaxed Diagram}
\label{sec:dd_refinement}
Simply having a relaxed DD is not enough to solve an optimization problem. In many problems, such as the TSP, the arc weights are exact, thus if the shortest path through the relaxed diagram is feasible, then it is also optimal. However, if the arc weights are just lower bounds, as in the ARP, the length of the shortest path might be super-optimal, even if the path is feasible.

\begin{figure}[!ht]
    \FIGURE{
        \centering
        \begin{subfigure}[b]{0.45\linewidth}
            \centering
             \vspace{10pt}
            \includegraphics[scale=.8]{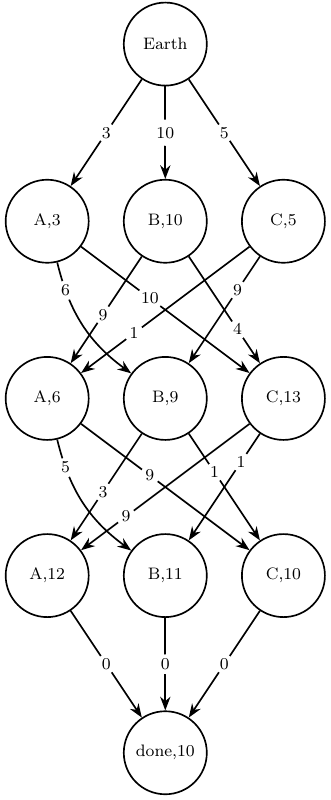}
            \captionsetup{font=normal}
            \caption{Initial Diagram}
            \label{fig:relaxed_dd_repeat1}
        \end{subfigure}
        \hfill
        \begin{subfigure}[b]{0.45\linewidth}
            \centering
             \vspace{10pt}
            \includegraphics[scale=.8]{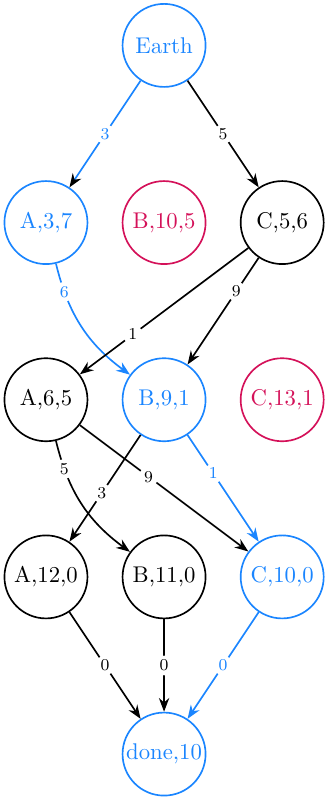}
            \captionsetup{font=normal}
            \caption{Constraint Propagation}
            \label{fig:relaxed_dd_3}
        \end{subfigure}
    }
    {Refining relaxed decision diagrams with constraint propagation.}
    {}
\end{figure}

\begin{figure}[!ht]
    \FIGURE{
        \centering
        \begin{subfigure}[b]{0.45\linewidth}
            \centering
             \vspace{10pt}
            \includegraphics[scale=.8]{Figures/WeightedInitialDD.pdf}
            \captionsetup{font=normal}
            \caption{Initial Diagram}
            \label{fig:relaxed_dd_repeat2}
        \end{subfigure}
        \hfill
        \begin{subfigure}[b]{0.45\linewidth}
            \centering
             \vspace{10pt}
            \includegraphics[scale=.8]{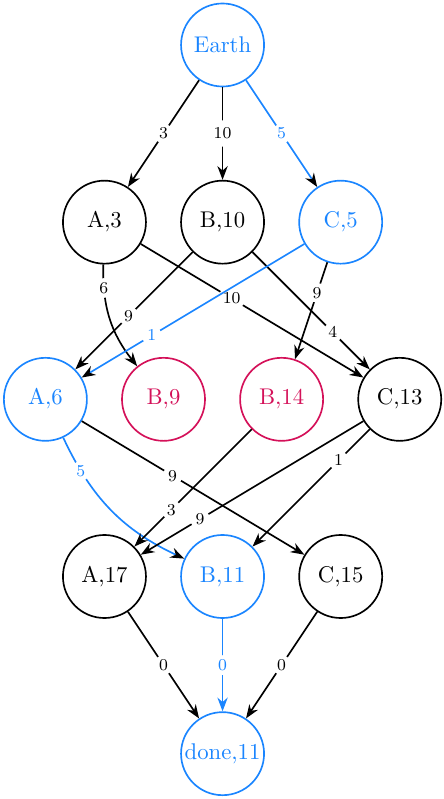}
            \captionsetup{font=normal}
            \caption{Node Splitting}
            \label{fig:relaxed_dd_4}
        \end{subfigure}
    }
    {Refining relaxed decision diagrams with node splitting.}
    {}
\end{figure}

\begin{figure}[!ht]
    \FIGURE{
        \centering
        \begin{subfigure}[b]{0.45\linewidth}
            \vspace{10pt}
            \centering
            \includegraphics[scale=.8]{Figures/WeightedInitialDD.pdf}
            \captionsetup{font=normal}
            \caption{Initial Diagram}
            \label{fig:relaxed_dd_repeat3}
        \end{subfigure}
        \hfill
        \begin{subfigure}[b]{0.5\linewidth}
             \vspace{10pt}
            \centering
            \includegraphics[scale=.8]{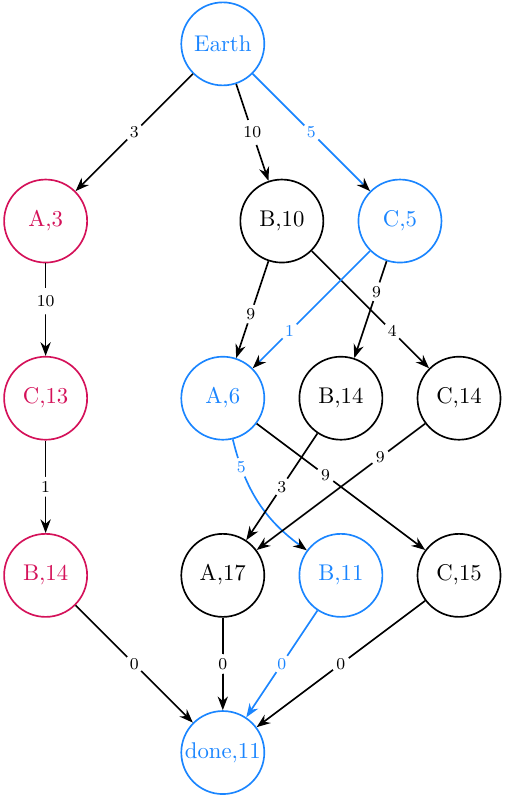}
            \captionsetup{font=normal}
            \caption{Node Peeling}
            \label{fig:relaxed_dd_5}
        \end{subfigure}
    }
    {Refining relaxed decision diagrams with node peeling.}
    {}
\end{figure}

To get an optimal solution to the ARP, we will leverage three DD refinement techniques. We will illustrate them using examples from before. Let $\mathcal{M}$ be our relaxed DD, let $\pi^*$ be the shortest path through $\mathcal{M}$, let $\underline{z}$ be the length of $\pi^*$ (and thus a relaxed bound on the optimal cost). In Figure \ref{fig:relaxed_dd_2}, we have a relaxed DD with a shortest path of $\pi^* = (A \rightarrow B \rightarrow C)$ and a cost of $\underline{z} = 10$. Assume that we solve the \textit{inner} problem of finding the optimal trajectories for the pairs in the sequence $\pi^*$, and find that $f(\pi^*) = 12$; now we know that our diagram is overestimating the quality of $\pi^*$, and we also have a real solution $\overline{z}$ with cost $\overline{z} = 12$. To make progress, we must now remove the path encoding $\pi^*$ from the DD without removing any different and potentially optimal solutions. This process is called refinement. In the remainder of this section, we discuss the refinement techniques and then demonstrate how they work using our example in Figures \ref{fig:relaxed_dd_3}, \ref{fig:relaxed_dd_4}, and \ref{fig:relaxed_dd_5}.

The first refinement technique we leverage is constraint propagation. When constraint propagation is used with DDs, it typically refers to adding constraints on the nodes or arcs such that if a constraint is violated, the node/arc can be removed. For example, each node $u$ in Figure \ref{fig:relaxed_dd_2} stores the length of the shortest path to that node $z_\downarrow(u)$. If each node similarly stores the length of the shortest path from $u$ to the terminal $z_\uparrow(u)$, then any node where $z_\downarrow(u)+z_\uparrow(u) \geq \overline{z}$ can be removed from the diagram. This is true because no path passes through such a node with a cost less than our best-known solution $\overline{z}$. \cite{CirHoe2013mdd} explores several such constraints for solving sequencing problems with DDs; our proposed algorithm includes all of them.

The second refinement technique we leverage is node splitting. In our example, we no longer need to include $\pi^*$ in our diagram because we know its true value. We must remove $\pi^*$ from $\mathcal{M}$ to find a tighter relaxed bound on our problem. This can be achieved with node splitting. Informally a node split: (1) creates a new node, (2) chooses a subset of the arcs pointing at the old node, and (3) redirects those arcs to point at the new node. 

Normally when a node $u$ is selected to be split, it is because a path being removed passes through $u$. Let $\phi$ be the label of the in-arc of $u$ on that path. To choose a subset of the arcs pointing at the old node $u$ during a split, we need a sorting function for the in-arcs of $u$. Let $All_v^\downarrow$ be the set of labels visited by every path from the root node $r$ to node $v$. For the ARP, an arc label is the destination asteroid represented by that arc in the transfer. Thus, $All_v^\downarrow$ is the set of asteroids visited on every path from $r$ to $v$. Let $a_{vu}$ be an arc from $u$ to $v$ and let $\arclabel(a)$ be the label on arc $a$. We use the function $h(a_{vu}, \phi)$ from \cite{CirHoe2013mdd} to determine which arcs to split on; it is defined as follows:
\begin{equation}
  \label{eq:sorting}
  h(a_{vu}, \phi) =\begin{cases}
    \textbf{true} & \text{if $\phi \in \{All_v^\downarrow\, \cup\, \arclabel(a_{vu})\}$}\\
    \textbf{false} & \text{otherwise}
  \end{cases}
\end{equation}
The result of using this function when splitting is that we will be more able to enforce the constraint that if an asteroid $\phi$ is visited on every path from the root node to the new node $u^\prime$, then $\phi$ cannot be revisited on any path from $u^\prime$ to the terminal. For example, $\phi$ may have been one of the labels on an out-arc of $u$, but any such out-arc of $u^\prime$ can be removed because it is infeasible. Node splitting is formalized in Algorithm \ref{algo:nodesplit}, and we explain our example in the next paragraph. 

\begin{algorithm}[!ht]
  \SetAlgoLined
  \DontPrintSemicolon
  \linespread{1.15}\selectfont
    \Input{a node $u$ on layer $L$ in relaxed DD $\mathcal{M}$}
    Let $in(u)$ be the the set of arcs that end at $u$\\
    Let $out(u)$ be the set of arcs that originate from $u$\\
    Let $a_{vu}$ be an arc going from node $v$ to node $u$\\
    Let $h(a_{vu},\phi)$ be a function that takes in an arc (as well as other relevant information such as $\phi$) and returns \textit{true} or \textit{false} (see Eq. \ref{eq:sorting})\;
    Create a new node $u^\prime$\;
    $L \assign L \cup \{u^\prime\}$\;

    \ForEach{$a_{vu} \in in(u)$}{
      \If(\tcp*[f]{For the ARP, $\phi = \nodelabel(u)$}){$h(a_{vu}, \phi)$}
            {Redirect $a$ such that $a_{v u^\prime}$}
    }
    \ForEach{$a_{uv}\in out(u)$}{
        Create and add arc $a_{u^\prime v}$ such that $\arclabel(a_{uv}) = \arclabel(a_{u^\prime v})$ \label{line:nodesplitarcs}\;
        filter($a_{u v}$), filter($a_{u^\prime v}$)\tcp*[r]{filter($a$) removes arcs that violate any constraint}
    }         
    \Return{$\mathcal{M}$}\;
    \caption{Node Splitting~\citep{CirHoe2013mdd}}
    \label{algo:nodesplit}
\end{algorithm}

A node $u$ is split by making a copy $u^\prime$. The in-arcs are heuristically sorted into two sets, and distributed between the two nodes. The out-arcs are copied so that both nodes have a copy of each arc. Then arcs are removed if they can be shown not to contain the optimal solution. In our example (Figure \ref{fig:relaxed_dd_repeat2}), node $(B,9)$ has two in-arcs $\{(A \rightarrow B), (C \rightarrow B)\}$, and two out-arcs $\{(B \rightarrow A), (B \rightarrow C)\}$. We split the node and give one in-arc to each, as well as a copy of both out-arcs. Then we observe that one node represents only paths starting with $(A \rightarrow B)$ and as explained at the beginning of Section \ref{sec:dd_refinement}, we want to remove the path $\pi^* = (A \rightarrow B \rightarrow C)$, so we can remove the $C$ out-arc for that node. We also observe that the other node represents only paths starting with $(C \rightarrow B)$, and the path $C \rightarrow B \rightarrow C$ is infeasible, so we can remove the $C$ arc from this node as well. Finally, the only path to the new $(B,9)$ node starts with $A$, and $A\rightarrow B \rightarrow A$ is infeasible, so we can remove the $A$ arc from this node, leaving it with no out-arcs. Any node with no out-arcs can simply be deleted because no potentially optimal paths pass through it.

The third refinement technique we leverage is the peel operation~\citep{RudCapRou2022cp,RudCapRou2023improved_pnb}. In the peel operation, we start by picking a node $u$, and then iteratively split nodes from the top down until we have separated all of the paths passing through $u$ in $\mathcal{M}$ into a discrete graph that only connects at the root and terminal. The process of repeatedly peeling and splitting nodes in a relaxed DD until you reach an optimal solution is called Peel-and-Bound (PnB). See Algorithm \ref{algo:Peel} for a formal treatment of the peel operation.

\begin{algorithm}[!ht]
\SetAlgoLined
  \DontPrintSemicolon
  \linespread{1.15}\selectfont
    \Input{a relaxed decision diagram $\mathcal{D}$ with root $r$ and terminal $t$, and a node $u$ in $\mathcal{D}$}
    \textbf{Output:} a relaxed decision diagram  $\mathcal{D}_u$ peeled from $\mathcal{D}$, and what remains of $\mathcal{D}$\\
    Let $in(u)$ for some node $u$ be the set of arcs that end at node $u$\\
    Let $out(u)$ be the set of arcs that originate from node $u$\\
    Let $in(\mathcal{D})$ be the set of arcs that end in diagram $\mathcal{D}$\\
    Let $out(\mathcal{D})$ be the set of arcs that originate in diagram $\mathcal{D}$\\
    Let $\mathcal{D}_u$ be an empty DD\;
    $in(u) \assign \emptyset$ \tcp*[l]{Remove the in-arcs of $u$}
    $\mathcal{D} \assign \mathcal{D}\backslash u$\;
    $\mathcal{D}_u \assign u$\;
    \While{$in(\mathcal{D})\cap out(\mathcal{D}_u) \neq \emptyset$}{
        \ForEach{node $m \in \mathcal{D}$ with an in arc that originates in $\mathcal{D}_u$}{
            Create a new node $m^\prime$ and add it to $\mathcal{D}_u$\;
            \ForEach{arc $a_{md} \in out(m)$}{
                Add arc $a_{m^\prime d}$\\
            }
            \ForEach{arc $a \in in(m)$ that originates in $\mathcal{D}_u$}{
                Change the destination of $a$ to $m^\prime$\;
                filter($a$)
            } 
            \ForEach{arc $a \in out(m)$}{
                filter($a$)
            }
        }
    }
    \While{$\exists m \in \mathcal{D}$ with $in(m)=\emptyset \lor out(m) = \emptyset$ (excluding $r$ and $t$)}{
            $in(m) \assign \emptyset$\;
            $out(m) \assign \emptyset$\;
            $\mathcal{D} \assign \mathcal{D} \backslash \{m\}$
    }
    \Return{($\mathcal{D}_u$, $\mathcal{D}$)}
    \caption{Peeling Procedure~\citep{RudCapRou2022cp,RudCapRou2023improved_pnb}}
    \label{algo:Peel}
\end{algorithm}

In our example (Figure \ref{fig:relaxed_dd_5}), we want to perform a peel operation that removes $\pi^* = (A \rightarrow B \rightarrow C)$ while also splitting the graph into a DD containing only paths that start with $A$, and a separate DD that contains no paths starting with $A$. To do this, we start by splitting the $B$ and $C$ node on the third layer into a version with $A$ as a parent, and a version that does not have $A$ as a parent. Then, the original $B$ node can be removed because it only contains $\pi^*$. Subsequently, the $B$ and $C$ nodes on the fourth layer are split so only paths that start with $A$ remain with the original nodes. Finally, the original $C$ nodes can be removed because no path starts with $A$ that does not also visit $C$. Because we remove $\pi^*$, the only remaining path starting with $A$ is $(A \rightarrow C \rightarrow B)$. Furthermore, the bound represented by the relaxed DD left behind, the one containing only nodes that start with $A$, is strengthened. Thus, the paths in the remaining diagram no longer intersect the ones in the peeled diagram.

\section{The Initial Decision Diagram}
\label{sec:initialconstruction}

This section details a technique for constructing a relaxed DD for the ARP that yields a valid relaxed bound.

\subsection{Initial Setup}
\label{sec:initialsetup}
We begin by constructing a relaxed DD $\mathcal{M}$ exactly as shown in Figure \ref{fig:relaxed_dd_1}. We start with a root node at layer $0$, a terminal node at layer $n+1$, and $n$ nodes on each layer with index $i$, where $1 \leq i \leq n$. Each node $u$ in a layer $i$ is labeled $\nodelabel(u)$ unique to that layer. Node labels are conceptually the same as arc labels. Each node and each arc has exactly one label, always representing either \Earth or an asteroid. An arc represents the transition from one node to another, so in the ARP, an arc represents a spacecraft transfer. The arc label is the asteroid being transferred to. A node represents a state and in the ARP, part of that state is the location of the spacecraft (either \Earth or an asteroid). Thus, the node label is either \Earth for the root or an asteroid for the other non-terminal nodes. The terminal node is a dummy node. Then, each node $v$ in layer $i-1$ is the origin for an arc ending at $u$, as long as $v$ has a different label. An arc with a null label is added from each node in layer $n$ to the terminal. The result is a relaxed DD that encodes every feasible sequence of asteroids as a path from the root to the terminal.

\subsection{Relaxing the Black-Box}
\label{sec:qrelaxed}

In real-world problems that depend on a black-box function or an inner problem, it is often possible to define relaxed versions of such functions that provide useful bounds or heuristic values \citep{EhmCamTho2016tdvrp,AndFagHob2015maritime}.

In the ARP, we can define a relaxed version of Eq.~\eqref{eq:inner} that removes the cost of waiting as follows:
\begin{equation}
  \label{eq:relaxedinner}
  \begin{split}
    f^\prime_\text{inner}(a, a', \eta, \tau, t) &= \Delta V + \frac{2\,\text{km/s}}{30\,\text{days}} \cdot t\\
    \text{s.t.}&\quad\tau \in [0, \tauf],\; t \in [1,\taumax]\\
    \text{where}&\quad\Delta V = \norm{\Delta\vec{v}_1} + \norm{\Delta\vec{v}_2}\quad\text{and}\quad (\Delta\vec{v}_1,\Delta\vec{v}_2) =
\text{Lambert}(a,a',\eta+\tau,t)
  \end{split}
\end{equation}

Since $f^\prime_\text{inner}$ removes the cost of waiting, we also relax the upper bound of $\tau$ to a parameter $\tauf$.
We can now define the following relaxed variant of $\mathcal{B}$:
\begin{equation}\label{eq:qprime}
  % \mathcal{B}(a, a', \eta) \to (\tau, t, z)
  \Bprime(a, a', \eta, \tauf)\colon A \times A \times \mathbb{R}^{+} \to [0, \tauf] \times [1, \tmax] \times \mathbb{R}^{+}
\end{equation}
which returns three values $(\tau, t, z)$, where $\tau$ and $t$ are the values that minimize $f^\prime_\text{inner}$ given $a$, $a'$, $\eta$ and $\tauf$, and $z$ is its minimal value. Note that this is only a relaxation in the sense that it removes part of the cost from the objective function, it does not change the method used to solve it. The procedure for solving $\Bprime$ is exactly the same as for solving $\mathcal{B}$.

Note that we will only ever use $\Bprime$ to compute bounds on costs, we will never use it to evaluate the cost of a permutation of asteroids. We always use $\mathcal{B}$ to evaluate the cost of solutions and partial solutions so that the values we report can be compared with the values reported by \cite{LopChiGil2022evo}.

\subsection{Calculating Valid Arc Bounds}
\label{sec:arcbounds}
The examples in Section \ref{sec:background} assumed the availability of bounds on the arc weights. In this section, we will explain how they can be calculated. We also mentioned that calculating arc bounds is expensive because it requires calculating the black-box function $\mathcal{B}$. In the ARP, this means that the inner problem must be solved using SLSQP.  The algorithm we propose is designed to minimize the times we need to calculate the value of $\mathcal{B}$.

When using constraint propagation, arcs are filtered and removed but never added. When splits and peels are performed, the number of arcs in the DD may grow, but only when the out-arcs of a node are copied. When this occurs, the asteroid that an arc is traveling from and the asteroid that the arc is traveling to remain unchanged. Therefore, when splitting or peeling, a valid bound on the cost of an arc is also a valid bound on its copy. While it may be desirable to recalculate the bound on an arc after a split to get a tighter bound, peels and splits can be carried out as many times as memory limitations allow without the need to re-evaluate $\mathcal{B}$. We propose a two-phase algorithm for creating an initial relaxed DD. This initial construction is intentionally computationally expensive because we will make enough evaluations of $\mathcal{B}$ that $\mathcal{B}$ does not need to be evaluated again after the initial relaxed DD is constructed, except when evaluating feasible solutions. 

\subsubsection{Phase One}\label{sec:phaseone}

Now we weight the DD. The arcs going from layer $0$ to layer $1$ are straightforward. A valid weight for an arc must be a valid lower bound on the cost of making the transfer represented by that arc. So a valid weight on arc $a_{uv}$ going from node $u$ to $v$ with labels $\nodelabel(u) = \Earth$ and $\nodelabel(v)$ is given by the $z$ value returned by $\mathcal{B}(\Earth, \nodelabel(v), \eta = 0)$. There is no possibility for $\eta$ other than $0$, so that $z$ value will be the smallest possible cost of going from Earth to $\nodelabel(v)$. 

In the ARP, these costs from layer 0 to layer 1 are exact because the inner problem is only solved for one pair of asteroids at a time for simplicity \citep{LopChiGil2022evo}. In other words, $\mathcal{B}$ does not depend on the path from node $v$ to the terminal. Thus, we can assume that the value of $\mathcal{B}$ for an arc only depends on the earliest start time at the starting node, which only depends on the path from the root node to that starting node. In the general trajectory optimization problem, it is possible to globally optimize the inner problem for a whole sequence of asteroids by optimizing all waiting and transfer times simultaneously. This simultaneous optimization may further reduce the cost of the sequence; for example, it may be that increasing the waiting time between early transfers leads to a reduction in the transfer cost of later transfers. In that case, the $\mathcal{B}$ value at a given arc may also depend on the path from the destination node to the terminal. 

Our proposed method works for either type of problem because the sum of the arc weights on a path is still a globally valid bound on the cost of the sequence that path represents. A dependence on the path from a node to the terminal simply makes the evaluation of feasible solutions more expensive because we need to re-optimize all waiting and transfer times simultaneously every time we weight an arc. In the ARP, if there is only one path from the root to a node $u$, and we know both the cost of the path from the root to $u$ and the earliest start time at the previous node, then the exact cost of that path can be calculated with a single evaluation of $\mathcal{B}$ that only optimizes the wait and transfer times of the arc that ends at $u$. In our implementation, we leverage this property by updating the weights of arcs on such a path with their exact cost and tracking the earliest start time for each node. This is not necessary to find an optimal solution, but it is a useful heuristic.

The arcs going from layers $i$ to $i+1$ for $i \in \{1,n-1\}$ remain to be weighted. We begin by finding a valid bound on the transfer between each ordered pair of asteroids $a,a' \in A$, with $a \neq a'$. We can use the no-wait $\Bprime$ (Eq.~\ref{eq:qprime}) for that purpose by noticing that the maximum time that any solution requires is trivially bounded by the sum of the upper bounds on the waiting times and travel times of the trajectories, i.e., $n(\taumax + \tmax)$. Therefore, the latest start time for the final transfer is bounded by $n(\taumax + \tmax) - \tmax$. That means that the transfer from an asteroid $a$ to an asteroid $a'$ can wait for $\tau \in [\earlieststart(a), n(\taumax + \tmax) - \tmax]$, where $\earlieststart(a)$ is the earliest start time at $a$. Initially, $\earlieststart(a) = \tau^{*} + t^{*}$, where $\tau^{*}$ and $t^{*}$ are the values returned by $\mathcal{B}(\Earth, a, \eta = 0)$ as described above. We can now define:
\newcommand{\zmin}{\ensuremath{z_\text{min}}}
\begin{equation}\label{eq:qmin}
  % \Bmin(a, a') := \Bprime (a, a', \tauf = n(\taumax + \tmax) - \tmax)
   \zmin(a, a') := z\quad\text{where}\quad (\tau, t, z) = \Bprime (a, a', \eta = \earlieststart(a), \tauf = n(\taumax + \tmax) - \tmax)
\end{equation}

Because $\Bprime$ does not penalize waiting time, the waiting time $\tau$ returned by $\Bprime$ above gives an optimal $\eta$ for minimizing the actual cost of the transfer using $\mathcal{B}$ (Eq.~\ref{eq:transfer}). In other words, the cost $\zmin(a,a')$  is the smallest possible cost of the transfer from $a$ to $a'$, and thus a lower bound of the $z$ value returned by $\mathcal{B}(a,a',\eta = \earlieststart(a) +\tau)$. Computing $\zmin(a_i,a_j)$ for all pairs $a_i,a_j$, where $i,j \in [1\dotsint n]$ and $i \neq j$ requires $n^2 - n$ calls to~$\Bprime$.

We add the pre-calculated $\zmin(a,a')$ values as arc weights. A valid weight for an arc $a_{uv}$ going from node $u$ with label $\nodelabel(u)$ to $v$ with label $\nodelabel(v)$, is simply $\zmin(\nodelabel(u),\nodelabel(v))$, which is the best possible transfer between $\nodelabel(u)$ and $\nodelabel(v)$ when ignoring waiting time at $\nodelabel(u)$. The resulting relaxed DD contains valid weights such that the sum of the costs of the arcs on any path from the root to the terminal provides a lower bound on the true cost of the permutation represented by the labels of the nodes on that path. However, this lower bound is quite weak, so we perform a second phase to strengthen the arc weights. 

\subsubsection{Phase Two}
\label{sec:phasetwo}
We will use a heuristic solution $\overline{z}$ to strengthen the arc weights in this phase. Many methods are available in the literature for finding a strong heuristic solution in trajectory optimization problems, such as the ARP \citep{SimIzzHaas2017multi,HenIzzLan2016fast}. In the ARP, a feasible solution is a valid permutation of the asteroids. For simplicity, our implementation uses an Euclidean distance heuristic to find the initial $\overline{z}$. It picks asteroids one by one, always moving to the closest unvisited asteroid as measured by Euclidean distance at the time of arrival to the last visited asteroid. This simple heuristic is exceptionally effective, and our proposed algorithm has no trouble rapidly finding strong feasible solutions.  In Section \ref{sec:heuristicsearch}, we will describe our method for finding feasible solutions in more detail.

Before we can use $\overline{z}$, we need to store some information on the nodes. We begin by recursively updating each node $u$ with the length of the shortest path from the root to $u$, denoted by $z_\downarrow(u)$, as well as the length of the shortest path from $u$ to the terminal, denoted by $z_\uparrow(u)$ \citep{CirHoe2013mdd}. Recall that the length of a path from the root to the terminal is a valid lower bound on the cost of the associated permutation. It follows that $z_\downarrow(u) + z_\uparrow(u)$ is a valid lower bound on the cost of any path that passes through $u$, because the shortest path that passes through $u$ is trivially the shortest path from the root to $u$ merged with the shortest path from $u$ to the terminal. Note that this implies that if $z_\downarrow(u) + z_\uparrow(u) > \overline{z}$, then the optimal solution cannot pass through $u$, and we can remove $u$ from the DD because the shortest path that passes through $u$ has a higher cost than a known feasible solution. 

Now consider an arc $a_{uv}$ connecting nodes $u$ to node $v$ with weight $w(a_{uv})$. Similarly, $z_\downarrow(u) + w(a_{uv}) +z_\uparrow(v)$ is a valid lower bound on the cost of any path that can pass through $a_{uv}$. This is because the shortest path passing through an arc is the shortest path to the origin of that arc, then the arc itself, then the shortest path from the destination of the arc to the terminal. As before, we know that if $z_\downarrow(u) + w(a_{uv}) + z_\uparrow(v) > \overline{z}$, then the optimal solution cannot pass through $a_{uv}$, and thus $a_{uv}$ can be removed from the graph. We will leverage this constraint to calculate new stronger arc weights.

Let us consider arc $a_{uv}$, and let the layer of $u$ be $1 \leq i \leq n-1$.
The latest possible start time $\eta$ for arcs going from layer $i$ to $i+1$ is $i \cdot (\taumax+\tmax)$. Thus, a naive approach to improving the quality of the weight of $a_{uv}$ is $\Bprime(\nodelabel(u), \nodelabel(v), \eta = \earlieststart(u), \tau_f = i \cdot (\taumax+\tmax) + \taumax)$, which returns the optimal transfer $(\tau, t, z)$ between two asteroids when waiting is free.

This $z$ value is a valid weight for each such arc $a_{uv}$, and greatly improves the weights on arcs on higher indexed layers over the weights calculated in phase one. However, this will do very little to improve the weights of arcs at lower indexed layers. This method is valid, but as we will see next, it is extremely inefficient. Consequently, we will now use the concepts we covered in this section to motivate a method that scales better with the layer index.

Recall the constraint that for arc $a_{uv}$, it must be that $z_\downarrow(u) + w(a_{uv}) + z_\uparrow(v) \leq \overline{z}$, or we can remove the arc. In the previous method, we used $i \cdot (\taumax+\tmax)$ as a bound on the free waiting time an arc can leverage. A better bound can be retrieved by recognizing that the true cost of $w(a_{uv})$ includes the cost of wait time as a component. Let $\tau_a$ be the wait time used by $a_{uv}$, and recall that the cost of time in the \emph{inner} objective function (Eq. \ref{eq:inner}) is $\frac{2}{30}$. It must also be the case that $z_\downarrow(u) + \frac{2}{30}\tau_a + z_\uparrow(v) \leq \overline{z}$. This is because if $\tau_a$ is large enough to cause that constraint to be violated, then $w(a_{uv})$ will also be large enough to cause its constraint to be violated because $w(a_{uv})$ is a sum of components that include $\frac{2}{30}\tau_a$. This yields the constraint that:
\begin{equation}
  \label{eq:waitcon}
  \tau_a \leq \frac{30}{2}(\overline{z} - z_\downarrow(u) - z_\uparrow(v))
\end{equation}
where the right-hand-side of Eq. \eqref{eq:waitcon} is a constant because $z_\downarrow(u)$ and $z_\uparrow(v)$ are known. This means that we can use this bound for $\tau_a$ to define a stronger bound per layer:
\begin{equation}\label{eq:strongbound}
\zmin(u, v, i) := z \:\:\: \text{where} \:\:\: (\tau, t, z) = \Bprime(\nodelabel(u), \nodelabel(v), \eta = \earlieststart(u) , \tauf) \:\:\: \text{and} \:\:\: \tauf = \tfrac{30}{2}(\overline{z} - z_\downarrow(u) - z_\uparrow(v))
\end{equation}

 Starting with $i = 1$, we calculate a new bound for every arc on a layer ($\zmin(u,v,i)$), then we update the values of $z_\downarrow(v)$ for each $v$ (in other words, each node on layer $i+1$), then we set $i \assign i+1$, and repeat until $i = n+1$. Therefore, we do not update the arcs with a weight of $0$ that go to the terminal. After this process is done, the $z_\uparrow$ values of the nodes need to be updated as they depend on the arcs' weights. 

Phase two requires performing a total of $n^3 - 2n^2 + n$ evaluations of $\Bprime$, but results in strong enough initial bounds that, after this initial setup, our algorithm for the ARP only needs to call the black-box function $\mathcal{B}$ to evaluate the true cost of partial feasible solutions; it does not need to evaluate $\mathcal{B}$ or $\Bprime$ to bound. 

The weights are improved based on the $z_\downarrow$ and $z_\uparrow$ values, but the values of $z_\downarrow$ and $z_\uparrow$ are updated based on the new arc weight. It is possible to run phase two repeatedly until the values converge. However, each additional update requires a huge amount of computation and, in practice, yields minuscule changes. So the decision of how many repetitions to perform while constructing the initial DD is heuristic. In our implementation, we found that this process rarely improves the bound enough to be useful, so we did not include it. 

\subsubsection{Earliest Start and Arrival Times}

So far, we have assumed that the black-box function, and thus the inner optimization, only depends on the transfer locations and the arrival time at the origin location (Eq.~\ref{eq:transfer}). In practice, this is often not the case, and it is possible to optimize all of the inner problem variables for a given permutation at once to find a better solution. More generally, the objective function $f(\pi)$ may be completely black-box and not linearly separable as in Eq.~\eqref{eq:objf}. Even in such cases, applying our proposed method may still be possible.

In the case of the ARP, and many global trajectory optimization problems, given a permutation of the asteroids or celestial bodies to be visited, it is possible in principle to globally optimize all times $\tau_i$ and $t_i$ for all $i=1,\dotsc,n$ in Eq.~\eqref{eq:objf}, and further improve the objective function value. The original definition of the ARP \citep{LopChiGil2022evo} does not perform this global inner optimization because it does not guarantee to find a better objective function value, making evaluating each permutation significantly more time-consuming. Nevertheless, in practical global trajectory optimization problems, performing such global inner optimization may be beneficial if enough computation time is available. Thus, we explain how to adapt our proposed method to handle such cases.

For the ARP as originally defined, calculating earliest start times is an unnecessary but useful heuristic. As explained in Section~\ref{sec:phaseone}, $\mathcal{B}(\Earth, \nodelabel(u), \eta =0)$ returns the $z$ value for an arc going from Earth to a node $u$ in layer $1$ as well as the exact $\tau+t$ value that the transfer requires, which is the exact earliest start time at $u$ for any valid sequence. However, if we optimize all of the inner problem decision variables for a valid sequence simultaneously, the earliest start time may be different. We describe next a valid alternative.

First, we define a variant of  $\mathcal{B}$ that limits the maximum total time $\theta$ for the transfer:
\begin{equation}
  \label{eq:qtilde}
  \begin{split}
    \Btilde(a, a', \eta, \theta) &= (\tau^*, t^*, z^*)\\
    \text{where}&\quad z^* = f_\text{inner}(a,a',\tau^*,t^*)\quad\text{and}\quad(\tau^*, t^*) = \argmin  f_\text{inner}(a,a',\tau,t)\\
    \text{s.t.}&\qquad\tau+t\leq \theta,\;\tau \in [0, \taumax],\; t \in [1,\tmax]\\
  \end{split}
\end{equation}

% \begin{equation}
%   \label{eq:qtilde}
%         \Btilde(a, a', \eta, \Theta_f)\colon A \times A \times [1, n(\taumax+\tmax) ] \to [0, \taumax] \times [1, \tmax] \times \mathbb{R}^{+}
% \end{equation}
% %
% \begin{equation}
%   \label{eq:qtilde}
%       \Btilde(a, a', \eta, \Theta_f)\colon A \times A \times [1, n(\taumax+\tmax) ] \to [0, \taumax] \times [1, \tmax] \times \mathbb{R}^{+}
% \end{equation}

The above function $\Btilde$ also returns three values $(\tau, t, z)$ with the same meaning as they have for $\mathcal{B}$, but the inner problem has the additional constraint that $\tau + t \leq \theta$.

Given a relaxed DD where the nodes store $z_\downarrow(u)$ and $z_\uparrow(u)$, we can use the above function to find the earliest start time ($\earlieststart$) for any node $u$ in layer 1. In other words, we can use $\Btilde$ to find the $\earlieststart$ for visiting an asteroid $\nodelabel(u)$ from Earth. This is accomplished by performing a binary search with a starting range of $[0,\taumax + \tmax]$ to find the value $\theta$ such that:
\begin{equation}\label{eq:slowsearch}
  \overline{z}-z_\uparrow(u) = \tilde{z}\quad\text{where}\quad(\tau, t, \tilde{z}) =  \Btilde(\Earth, a, \eta = 0, \theta)
\end{equation}

However, the $z$ values returned by $\mathcal{B}$, and thus $\Btilde$, are not necessarily smooth, so the equality may not be strictly satisfiable because the target value may not be in the range of \Btilde. When this is the case, the goal of the binary search should be modified to find the value of $\theta$ that produces the largest value of $\tilde{z} \leq \overline{z}-z_\uparrow(u)$.

The earliest arrival time ($\earliestarrival$) to reach the end of an arc $a_{uv}$ can be similarly found by performing a binary search to find the value of $\theta$ such that:
\begin{equation}
\tilde{z} = \overline{z}-z_\downarrow(u)-z_\uparrow(v)\quad\text{where}\quad(\tau, t, \tilde{z}) =
    \Btilde(\nodelabel(u), \nodelabel(v), \eta = \earlieststart(u), \theta) 
  \end{equation}

The $\earlieststart$ of a node $v$ is the minimum $\earliestarrival$ of arcs ending at $v$, which allows us to repeat the search for the nodes in the next layer. However, performing this search on a node $v$ requires running a binary search on every arc ending at $v$, which is significantly time-consuming. If using this method, we recommend only using it on nodes with one parent.

\section{Heuristic Search with Embedded Restricted Decision Diagrams}
\label{sec:heuristicsearch}
This section details a method for searching for feasible solutions to an ARP instance by reading the solutions from a relaxed DD.

\subsection{The General Case}
\label{sec:search_general_case}
To find solutions, we expand on the search procedure described by \cite{RudCapRou2023improved_pnb}, an extension of an idea from \cite{CopGilSch2024decision}. In this paper, we have discussed relaxed DDs, but have not included their much simpler mirror: restricted decision diagrams. A relaxed DD contains every permutation of decisions as a path from the root to the terminal, but accomplishes this by admitting infeasible permutations. A restricted DD is similarly constrained by a max width $\omega_s$, but it contains only feasible permutations of decisions as paths from the root to the terminal. This can be thought of as a generalized greedy heuristic, similar to beam search \citep{OwMor1988:ijpr}, that generates a fixed number of solutions and only continues to explore the ones with the lowest cost. Typical search procedures for DD-based solvers use relaxed DDs to generate lower bounds on minimization problems and separately construct restricted DDs to generate upper bounds. After generating a restricted DD, the shortest path from the root to the terminal $t$  is the best (lowest cost) solution found using the search procedure, and its length $z_\downarrow(t)$ is the solution cost.

A path in both a restricted DD and a relaxed DD represents a sequence of decisions. In a restricted DD, each path from the root to a node in the DD represents a feasible partial solution to the problem being solved. A relaxed DD contains every feasible solution, partial or otherwise. Thus, each path in a restricted DD will also exist in the associated relaxed DD. Furthermore, DDs are deterministic, so each path in a DD maps to exactly one node. Any possible path (sequence of decisions) in a restricted DD, will map to exactly one node in the associated relaxed DD. In other words, every possible restricted DD for an optimization problem will be embedded in a relaxed DD for that problem.

Formally, each path to a node $\overline{u}$ in a restricted DD  $\overline{\mathcal{M}}$ can be mapped to exactly one node $\underline{u}$ in a relaxed DD  $\underline{\mathcal{M}}$. Let the domain of a node $u$ be the set of feasible arc labels on out-arcs of $u$. When generating $\overline{\mathcal{M}}$, each node in $\overline{\mathcal{M}}$ creates a child node on the next layer for every element in its domain, and then the new layer is trimmed down to a pre-chosen maximum width.

The mapping from restricted DDs to relaxed DDs can be leveraged to improve the restricted DD as it is being generated. Let $d(\overline{u})$ be the domain of $\overline{u}$; set $d(\overline{u}) = d(\overline{u}) \cap d(\underline{u})$ before generating the child nodes of $\overline{u}$. This way, if an arc has been proven to be sub-optimal in $\underline{\mathcal{M}}$, it will not be created when generating $\overline{\mathcal{M}}$. 

This also allows for easy intensification of the search when combined with peel operations. If a node $u$ has been peeled from $\underline{\mathcal{M}}$ into $\underline{\mathcal{U}}$, then $\overline{\mathcal{M}}$ will not include any solutions that pass through $u$. Similarly, if a restricted DD $\overline{\mathcal{U}}$ is generated that is embedded in $\underline{\mathcal{U}}$, it will only explore solutions that pass through $u$. 

Without using this intersection operation as a way to trim the domain, $\overline{\mathcal{M}}$ will search the entire solution space that starts from the same root as $\underline{\mathcal{M}}$, even if the peeled diagrams have already been fully explored and are known to be sub-optimal. Using this method, each restricted DD will only search the solutions encoded within the matching relaxed DD. This means that each restricted DD has a significantly improved chance of finding the best solution embedded in the relaxed DD it maps to because the solution space it must explore becomes smaller with each peeled node. 

\subsection{Exploration Heuristic \& Usage for the ARP}
An additional benefit of leveraging embedded restricted DDs, and a critical benefit for solving ARPs, is that each node $u$ in a relaxed DD can also be labeled with the length of the shortest path from $u$ to the terminal, which we have already denoted as $z_\uparrow(u)$ in this paper. This provides a simple but powerful heuristic for deciding which nodes to explore. This heuristic is directly inspired by, but substantially different from, the rough bounding proposed in \cite{GilCopSch2021bbmdd}.

Let $\omega_s$ be the maximum width allotted to the search. Given a node to be processed $\overline{u}$ in a partially constructed restricted DD $\overline{\mathcal{M}}$, the typical method for figuring out which of $\overline{u}$'s children to explore would be to generate all of them and then repeatedly throw away the node $\overline{v}$ in the layer being constructed with the highest $z_\downarrow(v)$ until the width of the layer is $\omega_s$. However, for the ARP, each creation of a child requires a transfer from the last asteroid visited to an unvisited one, which requires a call to $\mathcal{B}$, causing such a search to become rapidly intractable.

For the embedded restricted DDs to be generated quickly, it is critical to limit the number of calls to $\mathcal{B}$. Evaluating a single feasible solution to the ARP requires $n-1$ calls to $\mathcal{B}$. We aim to perform a heuristic search that requires at most $\omega_s (n-1)$ calls to $\mathcal{B}$.

Let $\overline{V}$ be the set of potential child nodes that could be created during the procedure. Instead of generating every possible child node $\overline{v} \in \overline{V}$, we calculate a bound on the true cost of a solution passing through each $\overline{v}$. Then, we generate the $\omega_s$ child nodes in $\overline{V}$ with the lowest bounds. For a potential node $\overline{v}$ with parent $\overline{u}$, associated nodes in the relaxed DD $\underline{v}$ and $\underline{u}$ respectively, and arc $a_{\underline{u}\underline{v}}$ with weight $w(a_{\underline{u}\underline{v}})$, the bound $b(\overline{v})$ is calculated as:
\begin{equation} 
\label{eq:bound}
    b(\overline{v}) = z_\downarrow(\overline{u}) + w(a_{\underline{u}\underline{v}}) + z_\uparrow(\underline{v})
\end{equation}

For any node $\overline{u}$ in a restricted DD $\overline{\mathcal{M}}$, let $\underline{u}$ be the associated node in $\underline{\mathcal{M}}$. Let $\underline{r}$ be the root of an existing relaxed DD $\underline{\mathcal{M}}$. The procedure begins by creating $\overline{r}$ without creating its children. Then it iterates over the nodes in $\underline{\mathcal{M}}$ associated with the potential children of $\overline{r}$, and creates the $\omega_s$ best candidate children of $\overline{r}$, where candidacy is determined by the value returned by Eq. \ref{eq:bound}. Then this is repeated for each layer of the restricted DD. This procedure is formalized in Algorithm \ref{algo:search}.

\begin{algorithm}[bt]
  \caption{Heuristic Search Procedure}
  \label{algo:search}
\SetAlgoLined
\DontPrintSemicolon
    \Input{maximum width $\omega_s$; relaxed DD $\underline{\mathcal{M}}$ with root $r$ on layer $0$}
    Let $w(a_{xy})$ be the weight of the arc connecting $x$ to $y$.\;
    Let $\overline{\mathcal{M}}$ be a restricted DD with the same number of layers as $\underline{\mathcal{M}}$, where each layer is initialized to $\emptyset$ except layer $0$ which is initialized to $\{r\}$.\;
    For any node $\overline{u} \in \overline{\mathcal{M}}$, let $\underline{u}$ be the associated node in $\underline{\mathcal{M}}$.\;
    \ForEach{layer of nodes $L$ in $\overline{\mathcal{M}}$}{
        $q \assign \emptyset$ \tcp*[l]{$q$: empty list of nodes.}
        $b \assign \emptyset$ \tcp*[l]{$b$: empty mapping of nodes to values.}
        \ForEach{node $\overline{u} \in L$}{
            Let the domain of $\overline{u}$ be $d(\overline{u})$, initialized to be the set of node labels on the path from $r$ (Earth) to $\overline{u}$\;
            Let $d(\underline{u})$ be the set of node labels used by children of $\underline{u}$ \;
            $d(\overline{u}) \assign d(\overline{u}) \cap d(\underline{u})$\;
            \ForEach{label $l \in d(\overline{u}$)}{
                Create a new node $\overline{v}$ as a child of $\overline{u}$ with label $l$, but do not weight arc $a_{\overline{u}\overline{v}}$\;
                $b(\overline{v}) \assign z_\downarrow(\overline{u}) + w_{a_{\underline{u}\underline{v}}} + z_\uparrow(\underline{v})$\;
                $q \assign q \cup \{\overline{v}\}$\;
            }
        }
        Sort the elements $\overline{v} \in q$ from least to greatest by their bounds $b(\overline{v})$\;
        \While{$\abs{q} > \omega_s$}{
            Remove the last element of $q$ (the one with the largest bound) from both $q$ and $\overline{\mathcal{M}}$\;
        }
        \ForEach{node $\overline{v} \in q$}{
            $\tau^*, t^*, z^* \assign \mathcal{B}(l(\overline{u}), l(\overline{v}),\earlieststart(\overline{u}))$\;
            $w(a_{\overline{u}\overline{v}}) \assign z^*$\;
            $\earlieststart(\overline{v}) \assign \earlieststart(\overline{u}) +\tau^* + t^*$\;
        }
    }
    \Return{$\overline{\mathcal{M}}$}
\end{algorithm}

\section{Using Peel-and-Bound}
\label{sec:pnb}
Thus far in this paper, we have detailed a method of generating a valid initial relaxed DD for the ARP, strengthening the bounds on that DD with constraint propagation, splits, and peels, and searching that DD for improved feasible solutions. Those are all of the ingredients necessary to implement Peel-and-Bound. The process of putting them together is straightforward. After the insights presented in previous sections, we can adapt the algorithm from \cite{RudCapRou2023improved_pnb} for the ARP with the modification described below.

Let $\mathcal{M}(u)$ be a DD with root node $u$. After generating the initial relaxation $\underline{\mathcal{M}}(r)$ using the two-phase algorithm described in Section \ref{sec:initialconstruction}, place the entire DD into a queue $Q$ such that $Q = \{\underline{\mathcal{M}}(r)\}$. Then, a DD $\underline{\mathcal{M}}(u)$ is selected from $Q$ (for the first iteration $\underline{\mathcal{M}}(u) = \underline{\mathcal{M}}(r)$). Once a DD is selected, a search for feasible solutions is performed that respects the solution space defined by that DD (Algorithm \ref{algo:search}). Subsequently, a single exact node $e$ from $\underline{\mathcal{M}}(u)$ is selected. An exact node for the ARP is simply a node $u$ whose $z_\downarrow(u)$ is an exact value and not a bound. In the initial DD, all nodes on layer $1$ are exact because there is only one path from $r$ (Earth) to those nodes, and the true cost of those arcs is known. The process of selecting a DD and exact node are heuristic decisions discussed in Section \ref{sec:implementation}. 

Once $e$ has been selected, it is peeled from $\underline{\mathcal{M}}(u)$. We detailed this procedure in Algorithm \ref{algo:Peel}. Peeling $e$ accomplishes a top-down reading of the sub-graph induced by $e$, and potentially strengthens $\underline{\mathcal{M}}(u)$ by removing nodes and arcs in the process. If the shortest path through the modified $\underline{\mathcal{M}}(u)$ is less than the best-known solution, $\underline{\mathcal{M}}(u)$ is put back into $Q$. Then, the diagram $\underline{u}$ is strengthened using Algorithm \ref{algo:nodesplit}. Let $\underline{\mathcal{M}}(\underline{u})$ be the result; if the shortest path through the refined DD $\underline{\mathcal{M}}(\underline{u})$ is less than the best known solution, $\underline{\mathcal{M}}(\underline{u})$ is added to $Q$. The whole procedure is repeated until no DDs are left in the queue ($Q = \emptyset$). Peel-and-Bound is formalized in Algorithm \ref{algo:PnB}.

Peel-and-Bound will only terminate if no DDs are remaining that could possibly contain a solution better than the best-known solution. Therefore, when the algorithm terminates, the best-known solution is known to be optimal. However, it performs a search for feasible solutions at every iteration and may find (and return) the optimal solution long before that solution is known to be optimal. 

The only change required to use Peel-and-Bound for the ARP is in the decision of when to add a DD to the queue. In the original presentation of Peel-and-Bound, when the shortest path through a relaxed DD is a feasible solution, that DD is not placed back into the processing queue because the best-known feasible solution encoded in that DD is known. However, as explained in Section \ref{sec:dd_refinement}, the shortest path through a relaxed DD for the ARP being a feasible solution does not prove that path to be the optimal path in that DD because the arc weights are not exact. Thus, to refrain from adding a DD back into the queue for the ARP, there must be a known feasible solution with a cost that is not larger than the length of the shortest path in that DD. 

\begin{algorithm}[tb]
  \SetAlgoLined
  \DontPrintSemicolon
  \Input{The initial relaxed DD $\underline{\mathcal{M}}(r)$}
Let $v^*(u)$ be the shortest path that passes through $u$ (the lower bound encoded by only $u$).\;
 $Q \assign \{\underline{\mathcal{M}}(r)\}$ \;
 $\zopt \assign \infty$ \tcp*[l]{Value of the best known solution.}
 \While{$Q \neq \emptyset$}{
    $\mathcal{D} \assign$selectDiagram($Q$), $Q \assign Q\backslash\{\mathcal{D}\}$ \label{line:dd_select}\\
    $u \assign$selectExactNode($\mathcal{D}$)  \label{line:node_select}\;
    $\underline{u},\mathcal{D}^*\assign$ peel($\mathcal{D}$, $u$) \tcp*[l]{See Algorithm \ref{algo:Peel}}
    \If{$v^*(\mathcal{D}^*) < \zopt$}{
        $Q \assign Q \cup \{\mathcal{D}^*\}$
    }
    $\overline{\mathcal{M}} \assign \overline{\mathcal{M}}(u)$. \tcp*[l]{Using Algorithm \ref{algo:search}, with $u$ as the root}
    \If{$v^*(\overline{\mathcal{M}}) < \zopt$} {$\zopt \assign v^*(\overline{\mathcal{M}})$\;}
 	 \If{$\overline{\mathcal{M}}$ is not exact}{
 	 	$\underline{\mathcal{M}} \assign \underline{\mathcal{M}}(\underline{u})$\;
 	 	\If{$v^*(\underline{\mathcal{M}}) < \zopt$} {
 	 		$Q \assign Q \cup \{\underline{\mathcal{M}}\}$
 	 	}
 	 }
 }
 \Return{$\zopt$}
 \caption{Peel-and-Bound (PnB) Algorithm~\citep{RudCapRou2023improved_pnb}}
 \label{algo:PnB}
\end{algorithm}

\section{Implementation Details}
\label{sec:implementation}
This section details the data structures and heuristic decisions that we used to design an efficient implementation of the algorithms in this paper. 

\subsection[Memoization of B]{Memoization of $\mathcal{B}$}
The number of arcs stored in the processing queue during Peel-and-Bound can grow extremely rapidly. A core optimization when implementing Peel-and-Bound is refraining from performing unneeded writes to memory. As mentioned before, when new arcs are generated during Peel-and-Bound (line~\ref{line:nodesplitarcs} in Algorithm~\ref{algo:nodesplit}), they are copies of existing arcs and have the same cost. Therefore, we only want to store the cost once. Instead of storing costs on the arcs, we need a single location in memory where all such costs are stored. Storing feasible solutions can be accomplished straightforwardly, but storing arc bounds efficiently requires a more complex data structure. 

\subsubsection{Storing Feasible Solutions}
\label{sec:storingsols}

The arc costs for feasible solutions are stored in a basic decision tree to avoid redundant evaluations of $\mathcal{B}$. The tree is initialized with a root node $r$ having cost $c(r) = 0$ and a label of \Earth. When evaluating the cost of a partial (or complete) solution, the root node is queried for the first element (asteroid) $a$ in the sequence after \Earth. If the root node has a child $v$ whose label is also $a$, then the cost of going from Earth to $a$ is $c(v)$, and the next element in the sequence is compared with the children's labels of $v$. This process is repeated for each element of the given sequence, adding each individual node cost to a sum total until either the whole sequence is found in the tree. The total cost is returned or the next element $a$ in the sequence is not among the labels of the children of the corresponding tree node $u$. If $u$ does not have a child with the same label as $a$, then node $v$ is added to the tree and its fields are set as follows:
\begin{equation}
\begin{split}
\tau^*, t^*, z^* &\assign \mathcal{B}(\nodelabel(u)), \nodelabel(v),\earlieststart(u))\\
            c(v) &\assign z^*\\
            \earlieststart(v)&\assign \earlieststart(u) +\tau^* + t^*
\end{split}
\end{equation}

The above process ensures that the cost of any feasible solution, partial or complete, is only computed once. Although a complete tree for an instance of size $n$ has $1 + \sum_{i=1}^n\prod_{j=0}^{i-1} (n-j)$ nodes, in practice, only a very small part of the tree is generated. Even with a complete tree where the children of each node are kept sorted, it would take at most $\mathcal{O}(k \log n)$ comparisons to retrieve the cost of a sequence of length $k \leq n$, which is much faster than evaluating the sequence.

\subsubsection{Memorizing Arc Bounds}
\label{sec:storingbounds}
During Peel-and-Bound, we need to calculate arc weights every time we peel a DD (Algorithm \ref{algo:Peel}), and every time we need to compute the shortest path through a DD. However, we would like to reuse the values already computed by the  $\zmin$ functions of phases one and two (Eqs.~\ref{eq:qmin} and~\ref{eq:strongbound}, respectively). Thus, we wish to look up a pre-computed value of $\zmin$ that would be valid for the new arc. 

Both definitions of $\zmin$ only depend on the evaluation of $\Bprime(a, a', \eta, \tauf)$. Moreover, by the definition of the inner optimization problem (Eq.~\ref{eq:relaxedinner}), we also know that $\Bprime(a, a', \eta', \tauf')$ cannot return a cost $z$  better than $\Bprime(a, a', \eta, \tauf)$ if $[\eta', \eta'+\tauf'] \subseteq [\eta, \eta+\tauf]$, because the inner optimization problem is the same but with tighter bounds on $\tau$ for the former.

During a run of Peel-and-Bound, the time intervals $[\eta, \eta+\tauf]$ implied by the arcs can only decrease in the above manner. Let $u$ be a node split into $u^\prime_1$ and $u^\prime_2$, and let $c(u)$ be the cost of node $u$. The bounds on the new nodes must be at least as tight as the bounds on the old nodes, so $c(u) \leq c(u^\prime_1), c(u) \leq c(u^\prime_2)$. The earliest starting time of the nodes can increase but not decrease, so $\earlieststart(u) \leq \earlieststart(u^\prime_1)$, $\earlieststart(u) \leq \earlieststart(u^\prime_2)$. This means that the interval start can get larger but not smaller. Let $\tau_u$ be the wait time associated with $u$, the wait time of the peeled nodes can decrease but not increase, $\tau_u \geq \tau_{u^\prime_1}$, $\tau_u \geq \tau_{u^\prime_2}$, which means that the interval end can get smaller but not larger. Thus, the time interval $[\eta, \eta+\tauf]$ used to compute the arc weight associated with a node $u$, must contain the time intervals that would be required to compute the arc weights associated with $u^\prime_1$ and $u^\prime_2$.

Let us store in an interval tree the $\zmin$ values resulting from any evaluation of $\Bprime$ with origin $a$ and destination $a'$ together with the time intervals used in the evaluation. Given an arc with an associated time interval and whose origin and destination correspond to $a$ and $a'$, every stored $\zmin$ value with a time interval containing the arc's time interval is a valid bound for the arc. The highest value among those stored provides the strongest bound.

Interval trees are a standard data structure designed for efficiently finding all intervals that overlap with any given interval or point. They have a lookup time for finding overlapping intervals in $\mathcal{O}(n+m)$ where $n$ is the number of intervals in the tree, and $m$ is the number of intervals that satisfy the query. A standard interval tree is optimized for finding overlapping intervals; each tree-node stores the maximum value of its descendants. In our case, we only need to find containing intervals. Thus, we also store its descendants' minimum value (interval start) in each tree node. This prevents wasted time from looking for containing intervals on a branch of the tree when the lower bound on that branch is higher than the start time of our interval. We also adjusted the lookup function to only return containing intervals instead of all overlapping intervals.

\subsection{Heuristic Decisions}\label{sec:heuristics}
In our implementation, we test two different methods for selecting a DD from $Q$ (line \ref{line:dd_select} in Algorithm \ref{algo:PnB}): worst bound and largest index. Selecting the DD with the worst bound is similar to performing a breadth-first-search; it is deciding to always improve the current lower bound on the problem but can lead to a very large $Q$. Selecting the DD with a root node whose index in the initial DD was highest is similar to performing a depth-first search; it prioritizes keeping $Q$ small over improving the lower bound as quickly as possible. 

We also test two different methods for selecting a node to peel from DD $Q$ (line \ref{line:dd_select} in Algorithm \ref{algo:PnB}). The first method is to pick the last exact node, the highest indexed exact node, on the shortest path through $Q$. The second method is to pick the maximal node, the lowest indexed node on the shortest path through $Q$. These methods are discussed in detail in \cite{RudCapRou2023improved_pnb}.

\subsection{Limitations of the Inner Optimizer}
\label{sec:limits}

The inner problem is non-convex and non-linear, making it impossible, in general, to prove that a local optimum is also a global optimum. In $\mathcal{B}$, the search space is small enough that local optima are also likely to be global optima. Even if this is not the case, the inner optimizer is deterministic, meaning that the same permutation will yield the same local optima for the inner optimization problems. The goal is to find the globally optimal permutation. When using $\mathcal{B'}$, the inner optimizer searches within a superset of the search space induced by $\mathcal{B}$, potentially finding a local optimum that is worse than the one returned by $\mathcal{B}$. In such cases, the outer optimization might be misled about which permutation is optimal. However, the permutation found is still a solution (possibly suboptimal) to the original ARP, because feasible solutions are always evaluated using $\mathcal{B}$, not $\mathcal{B'}$.

More concretely, the inner optimization that computes $\mathcal{B}$ in \cite{LopChiGil2022evo} consists of a single deterministic run of SLSQP, as explained in Section~\ref{sec:ARP}. A single run of SLSQP almost always returns a global optimum when $\tauf \leq \taumax$  as set in the original paper. However, SLSQP may return a locally optimal solution when $\tauf \gg \taumax$,  as we set it here when computing $\Bprime$ (Eq.~\ref{eq:qprime}). That is, $\Bprime(a,a',\eta = 0,\tauf)$ may return $(\tau',t', z')$ as optimal yet $\Bprime(a,a', \eta=\tau', \tauf)$ may return a better solution even though the time interval $[\tau', \tauf]$ is contained within $[0, \tauf]$. This problem can be alleviated by performing multiple restarts of SLSQP for each evaluation of $\Bprime$. In other words, when optimizing the relaxed inner problem (Eq.~\ref{eq:relaxedinner}), we can partition the range $[\eta, \eta+\tauf]$ into disjoint intervals, run SLSQP for each interval, and keep the best solution. We report experiments with different numbers of restarts that show their effect on runtime and solution quality. 

SLSQP could be replaced with a more effective optimizer, such L-BFGS-B \citep{NocWri2006} or CMA-ES \citep{HanOst2001ec}. For other inner problems, it may be possible to ensure that the inner optimizer returns a global optimum. In this paper, we decided to keep SLSQP as the inner optimizer to compare our results on the ARP with those reported by \cite{LopChiGil2022evo}. The following experiments (Section \ref{sec:exps}) show that increasing the quality of the inner optimizer does lead to better solutions, but not always, and not by much, so it is likely that we are already finding the optimal solutions for several of the ARP instances. 

\section{Experimental Results}
\label{sec:exps}

An ARP instance with $n$ asteroids is constructed by randomly selecting $n$ asteroids from a database containing $83,453$
asteroids. The selection process is controlled by a \emph{seed} value, ensuring that specific instances can be consistently replicated. In \cite{LopChiGil2022evo}, they test their algorithms with $n \in \{10,15,20,25,30\}$ and \emph{seed} $\in \{42,73\}$. We ran experiments with the same $n$ values but added three randomly chosen seeds: $8$, $22$, and $59$. The experiments were performed on a computer with an AMD Rome 7532 at 2.40 GHz and 64Gb RAM. Our code and raw experimental data are available at \url{https://zenodo.org/doi/10.5281/zenodo.12675217}. Several results tables in this section have a \emph{queue} column. This value is the number of DDs remaining in the processing queue when the solver ran out of time.

\newcommand{\pmulti}{\textit{multi}}

\subsection{Note on Optimality}
As discussed in Section \ref{sec:limits}, changing the quality of the inner optimizer may impact the quality of solutions returned by our framework. In this section, we report different optimal solutions using different settings for several problems. These different solutions result from the inner optimizer not being of sufficient quality for some settings, and thus are intentional and do not represent a mistake in the implementation. \pmulti \: is the setting we used to change the quality of the inner solver. It is the number of restarts that the inner solver uses when searching for an optimal solution.

\subsection{Initial Experiment: Determining Best Settings}
In our first test, shown in Table \ref{tab:test1_15}, we sought to determine which peel setting would be the most effective, which DD-widths $\omega$ would be the most effective, and the effect of using depth-first search (dfs, see Section~\ref{sec:heuristics}). The value of $\omega$ used can greatly impact the solve time because relaxed DDs with larger widths generally yield better bounds but also take longer to compute. The best value for $\omega$ is the one that balances this quality/time trade-off. 

\begin{figure}[tb]
  \centering%
  \includegraphics[width=\textwidth]{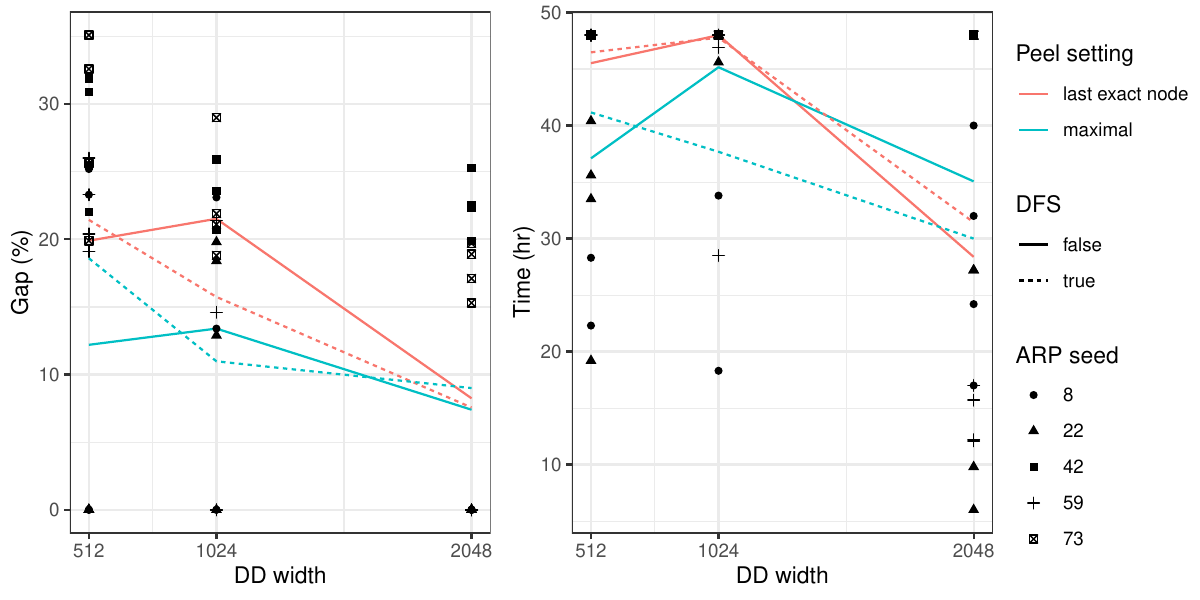}
  \caption{ARP instances with $n=15$, a 2 day runtime, an embedded search width of 100, and $\pmulti=1$ (Table~\ref{tab:test1_15})}
    \label{fig:test1_15}
  \end{figure}

We tested $n=15$ with all 5 seeds, 2 days of runtime, and no SLSQP restarts ($\pmulti = 1$ means one SLSQP run for solving the inner optimization problem). Figure~\ref{fig:test1_15} summarizes these results; the lines show the mean gap (or mean time) over the ARP instances.

\begin{table}[!p]
\caption{ARP instances with $n=15$, a 2 day runtime, an embedded search width of 100, and $\pmulti=1$}\label{tab:test1_15}
\centering%
\resizebox*{!}{0.9\textheight}{%
\begin{tabular}{|c|ccc|ccc|c|}
\toprule
seed                 & peel setting                     & dfs                    & DD width & lb    & ub    & gap (\%) & time (hr) \\\midrule
\multirow{12}{*}{8}  & \multirow{3}{*}{last exact node} & \multirow{3}{*}{true}  & 512      & 406.0 & 528.9 & 23.3     & -         \\
                     &                                  &                        & 1024     & 433.3 & 500.4 & 13.4     & -         \\
                     &                                  &                        & 2048     & 499.5 & 499.5 & 0.0      & 17.0      \\\cline{2-8} 
                     & \multirow{3}{*}{maximal}         & \multirow{3}{*}{true}  & 512      & 499.5 & 499.5 & 0.0      & 28.3      \\
                     &                                  &                        & 1024     & 500.1 & 500.1 & 0.0      & 18.3      \\
                     &                                  &                        & 2048     & 504.5 & 504.5 & 0.0      & 32.0      \\\cline{2-8} 
                     & \multirow{3}{*}{last exact node} & \multirow{3}{*}{false} & 512      & 411.3 & 550.0 & 25.2     & -         \\
                     &                                  &                        & 1024     & 412.6 & 536.5 & 23.1     & -         \\
                     &                                  &                        & 2048     & 500.9 & 500.9 & 0.0      & 24.2      \\\cline{2-8} 
                     & \multirow{3}{*}{maximal}         & \multirow{3}{*}{false} & 512      & 499.5 & 499.5 & 0.0      & 22.3      \\
                     &                                  &                        & 1024     & 499.5 & 499.5 & 0.0      & 33.8      \\
                     &                                  &                        & 2048     & 503.8 & 503.8 & 0.0      & 40.0      \\\midrule
\multirow{12}{*}{22} & \multirow{3}{*}{last exact node} & \multirow{3}{*}{true}  & 512      & 472.4 & 472.4 & 0.0      & 40.4      \\
                     &                                  &                        & 1024     & 400.9 & 500.0 & 19.8     & -         \\
                     &                                  &                        & 2048     & 501.7 & 501.7 & 0.0      & 27.2      \\\cline{2-8} 
                     & \multirow{3}{*}{maximal}         & \multirow{3}{*}{true}  & 512      & 472.4 & 472.4 & 0.0      & 33.5      \\
                     &                                  &                        & 1024     & 501.7 & 501.7 & 0.0      & 45.6      \\
                     &                                  &                        & 2048     & 481.9 & 481.9 & 0.0      & 9.8       \\\cline{2-8} 
                     & \multirow{3}{*}{last exact node} & \multirow{3}{*}{false} & 512      & 472.4 & 472.4 & 0.0      & 35.6      \\
                     &                                  &                        & 1024     & 406.6 & 498.5 & 18.4     & -         \\
                     &                                  &                        & 2048     & 482.3 & 482.3 & 0.0      & 6.0       \\\cline{2-8} 
                     & \multirow{3}{*}{maximal}         & \multirow{3}{*}{false} & 512      & 472.4 & 472.4 & 0.0      & 19.2      \\
                     &                                  &                        & 1024     & 438.9 & 504.0 & 12.9     & -         \\
                     &                                  &                        & 2048     & 501.7 & 501.7 & 0.0      & 27.2      \\\midrule
\multirow{12}{*}{42} & \multirow{3}{*}{last exact node} & \multirow{3}{*}{true}  & 512      & 351.2 & 508.2 & 30.9     & -         \\
                     &                                  &                        & 1024     & 382.3 & 500.1 & 23.6     & -         \\
                     &                                  &                        & 2048     & 393.4 & 507.5 & 22.5     & -         \\\cline{2-8} 
                     & \multirow{3}{*}{maximal}         & \multirow{3}{*}{true}  & 512      & 346.0 & 508.2 & 31.9     & -         \\
                     &                                  &                        & 1024     & 376.7 & 508.2 & 25.9     & -         \\
                     &                                  &                        & 2048     & 379.6 & 508.2 & 25.3     & -         \\\cline{2-8} 
                     & \multirow{3}{*}{last exact node} & \multirow{3}{*}{false} & 512      & 379.3 & 508.2 & 25.4     & -         \\
                     &                                  &                        & 1024     & 8.1   & 508.2 & 23.6     & -         \\
                     &                                  &                        & 2048     & 394.5 & 507.5 & 22.3     & -         \\\cline{2-8} 
                     & \multirow{3}{*}{maximal}         & \multirow{3}{*}{false} & 512      & 395.3 & 506.9 & 22.0     & -         \\
                     &                                  &                        & 1024     & 402.2 & 506.9 & 20.7     & -         \\
                     &                                  &                        & 2048     & 406.0 & 506.9 & 19.9     & -         \\\midrule
\multirow{12}{*}{59} & \multirow{3}{*}{last exact node} & \multirow{3}{*}{true}  & 512      & 416.6 & 523.2 & 20.4     & -         \\
                     &                                  &                        & 1024     & 523.2 & 523.2 & 0.0      & 46.9      \\
                     &                                  &                        & 2048     & 530.4 & 530.4 & 0.0      & 17.0      \\\cline{2-8} 
                     & \multirow{3}{*}{maximal}         & \multirow{3}{*}{true}  & 512      & 403.7 & 545.8 & 26.0     & -         \\
                     &                                  &                        & 1024     & 523.2 & 523.2 & 0.0      & 28.5      \\
                     &                                  &                        & 2048     & 530.4 & 530.4 & 0.0      & 12.2      \\\cline{2-8} 
                     & \multirow{3}{*}{last exact node} & \multirow{3}{*}{false} & 512      & 420.9 & 548.9 & 23.3     & -         \\
                     &                                  &                        & 1024     & 426.3 & 542.7 & 21.4     & -         \\
                     &                                  &                        & 2048     & 530.4 & 530.4 & 0.0      & 15.7      \\\cline{2-8} 
                     & \multirow{3}{*}{maximal}         & \multirow{3}{*}{false} & 512      & 440.3 & 543.9 & 19.1     & -         \\
                     &                                  &                        & 1024     & 446.7 & 523.2 & 14.6     & -         \\
                     &                                  &                        & 2048     & 530.4 & 530.4 & 0.0      & 12.1      \\\midrule
\multirow{12}{*}{73} & \multirow{3}{*}{last exact node} & \multirow{3}{*}{true}  & 512      & 333.9 & 495.2 & 32.6     & -         \\
                     &                                  &                        & 1024     & 389.2 & 498.0 & 21.9     & -         \\
                     &                                  &                        & 2048     & 419.5 & 495.2 & 15.3     & -         \\\cline{2-8} 
                     & \multirow{3}{*}{maximal}         & \multirow{3}{*}{true}  & 512      & 326.3 & 502.5 & 35.1     & -         \\
                     &                                  &                        & 1024     & 351.7 & 495.2 & 29.0     & -         \\
                     &                                  &                        & 2048     & 402.3 & 500.8 & 19.7     & -         \\\cline{2-8} 
                     & \multirow{3}{*}{last exact node} & \multirow{3}{*}{false} & 512      & 373.8 & 502.5 & 25.6     & -         \\
                     &                                  &                        & 1024     & 396.3 & 502.5 & 21.1     & -         \\
                     &                                  &                        & 2048     & 407.6 & 502.5 & 18.9     & -         \\\cline{2-8} 
                     & \multirow{3}{*}{maximal}         & \multirow{3}{*}{false} & 512      & 398.7 & 498.1 & 19.9     & -         \\
                     &                                  &                        & 1024     & 408.1 & 502.5 & 18.8     & -         \\
                     &                                  &                        & 2048     & 416.3 & 502.1 & 17.1     & -         \\\bottomrule
\end{tabular}%
}
\\[10pt] % Adds some space after the table
{\small
    \emph{dfs} = depth first search (queue processing order),
    \emph{DD width} = width of relaxed DD, \\\emph{lb} = lower bound, \emph{ub} = upper bound, \emph{gap} = optimality gap, \emph{time} = time to solve
}
\end{table}

We tested two peel settings: maximal and last exact node. Although there is no clear winner, maximal generally performs better and performs only slightly worse when it does not. Consequently, we will use the maximal setting in all future experiments. The effectiveness of the depth-first-search (dfs) queue ordering appears random based on summary data. The raw data, which is also available in the repository, contains a timestamp and bounds for each improvement to the bounds made during a run of PnB. The raw data suggests the success of dfs mainly depends on whether the best solution is found in an early branch. It also suggests that dfs is slightly worse for the \emph{maximal} peel setting and both $\omega=512$ and $\omega=2048$, as also shown in the left plot of Fig.~\ref{fig:test1_15}. We disabled dfs ordering in the remaining experiments, opting for the default ordering that processes the DD with the weakest bound, likely reducing the optimality gap in unresolved problems. Regarding relaxed DD-width, we tested 512, 1024, and 2048. While 2048 showed the greatest success in closing problems and narrowing the optimality gap, 512 resolved specific settings significantly faster.  Moving forward, we will retest 512 and 2048, and introduce a new width of 256. We do not include tests with widths larger than 2048 because in our experience with DDs, the slowdown in construction time from increasing the width starts to become substantial around 2048. This happens because the number of arcs at each layer is at most $\mathcal{O}(\omega^2 n)$, and small increases in the width beyond 2048 can lead to enormous increases in the construction time.

\begin{figure}[tb]
  \centering%
  \includegraphics[width=\textwidth]{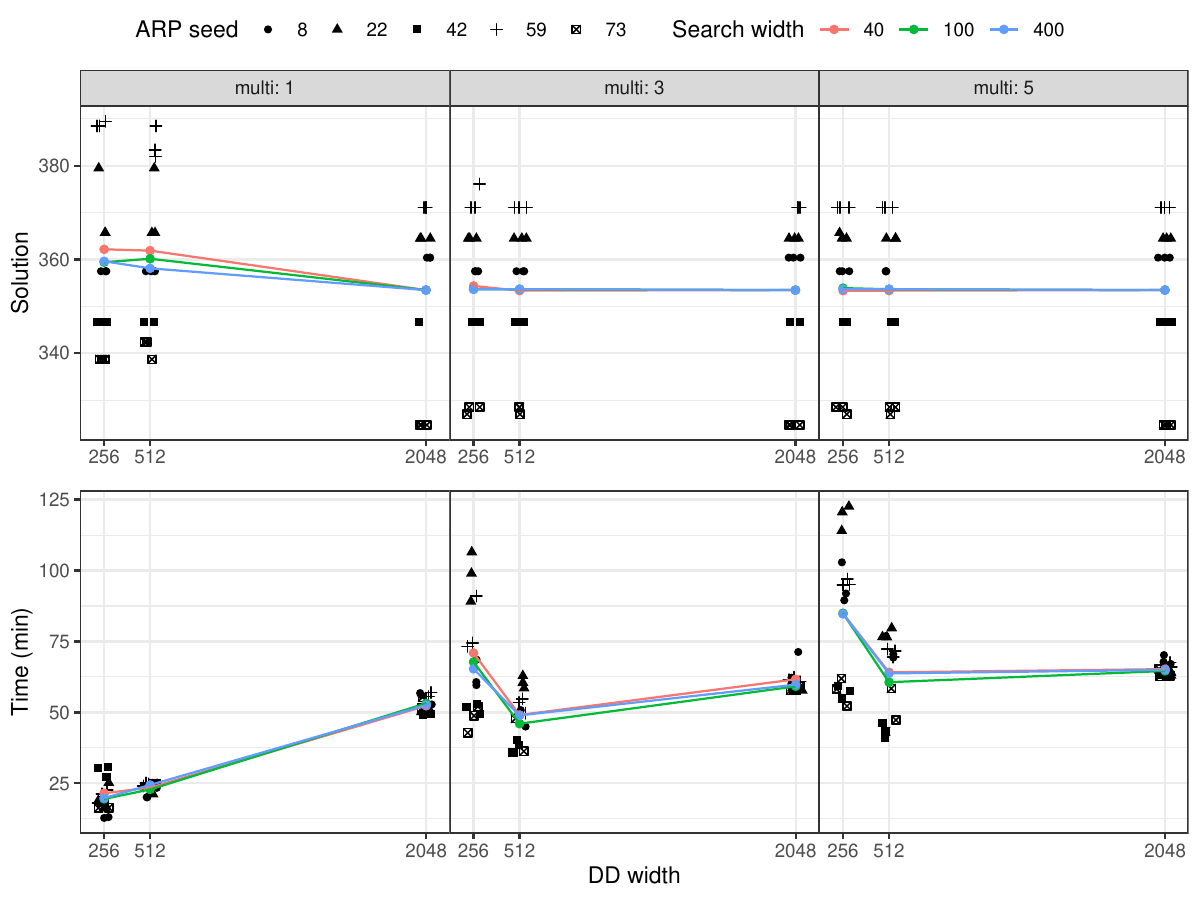}
  \caption{ARP instances with $n=10$. Lines show mean value over ARP instances (Table~\ref{tab:test2_10})}
  \label{fig:test2_10}
\end{figure}

\begin{table}[tb]
\caption{ARP instances with $n=10$}\label{tab:test2_10}
\centering
\resizebox{!}{6.75cm}{%
\begin{tabular}{|c|rr|cc|cc|cc|}
\toprule
\multirow{2}{*}{seed} & \multirow{2}{*}{search} & \multirow{2}{*}{DD width} & \multicolumn{2}{c|}{$\pmulti=1$} & \multicolumn{2}{c|}{$\pmulti=3$} & \multicolumn{2}{c|}{$\pmulti=5$}              \\
                      &                         &                           & solution                       & time(min)                      & solution & time(min) & solution & time(min) \\\midrule
\multirow{9}{*}{8}    & \multirow{3}{*}{40}     & 256                       & 357.5                          & 13.1                           & 357.5    & 60.7      & 357.5    & 89.5      \\
                      &                         & 512                       & 357.5                          & 20.1                           & 357.5    & 49.6      & 357.5    & 69.3      \\
                      &                         & 2048                     & 360.4                          & 52.8                           & 360.4    & 71.3      & 360.4    & 70.2      \\\cline{2-9} 
                      & \multirow{3}{*}{100}    & 256                       & 357.5                          & 12.8                           & 357.5    & 59.6      & 357.5    & 91.9      \\
                      &                         & 512                       & 357.5                          & 20.1                           & 357.5    & 45.0      & 357.5    & 70.5      \\
                      &                         & 2048                     & 360.4                          & 56.1                           & 360.4    & 60.0      & 360.4    & 67.8      \\\cline{2-9} 
                      & \multirow{3}{*}{400}    & 256                       & 357.5                          & 16.1                           & 357.5    & 68.6      & 357.5    & 102.9     \\
                      &                         & 512                       & 357.5                          & 23.4                           & 357.5    & 50.8      & 357.5    & 76.9      \\
                      &                         & 2048                     & 360.4                          & 56.8                           & 360.4    & 62.2      & 360.4    & 67.1      \\\midrule
\multirow{9}{*}{22}   & \multirow{3}{*}{40}     & 256                       & 379.5                          & 25.2                           & 364.5    & 99.0      & 364.5    & 120.6     \\
                      &                         & 512                       & 379.5                          & 24.5                           & 364.5    & 58.5      & 364.5    & 79.7      \\
                      &                         & 2048                     & 364.5                          & 50.2                           & 364.5    & 57.7      & 364.5    & 63.0      \\\cline{2-9} 
                      & \multirow{3}{*}{100}    & 256                       & 365.7                          & 16.2                           & 364.5    & 106.5     & 365.7    & 122.6     \\
                      &                         & 512                       & 365.7                          & 21.2                           & 364.5    & 60.5      & 364.5    & 76.6      \\
                      &                         & 2048                     & 364.5                          & 51.4                           & 364.5    & 58.5      & 364.5    & 64.1      \\\cline{2-9} 
                      & \multirow{3}{*}{400}    & 256                       & 365.7                          & 18.9                           & 364.5    & 89.1      & 364.5    & 114.0     \\
                      &                         & 512                       & 365.7                          & 23.8                           & 364.5    & 62.9      & 364.5    & 76.5      \\
                      &                         & 2048                     & 364.5                          & 51.2                           & 364.5    & 58.8      & 364.5    & 64.0      \\\midrule
\multirow{9}{*}{42}   & \multirow{3}{*}{40}     & 256                       & 346.7                          & 30.3                           & 346.7    & 51.9      & 346.7    & 57.5      \\
                      &                         & 512                       & 346.7                          & 23.7                           & 346.7    & 36.0      & 346.7    & 40.9      \\
                      &                         & 2048                     & 346.7                          & 50.1                           & 346.7    & 59.3      & 346.7    & 63.9      \\\cline{2-9} 
                      & \multirow{3}{*}{100}    & 256                       & 346.7                          & 30.9                           & 346.7    & 49.5      & 346.7    & 54.9      \\
                      &                         & 512                       & 346.7                          & 24.8                           & 346.7    & 38.6      & 346.7    & 43.3      \\
                      &                         & 2048                     & 346.7                          & 49.3                           & 346.7    & 58.0      & 346.7    & 62.4      \\\cline{2-9} 
                      & \multirow{3}{*}{400}    & 256                       & 346.7                          & 27.2                           & 346.7    & 52.9      & 346.7    & 59.3      \\
                      &                         & 512                       & 346.7                          & 25.1                           & 346.7    & 40.4      & 346.7    & 46.3      \\
                      &                         & 2048                     & 346.7                          & 49.1                           & 346.7    & 57.7      & 346.7    & 63.3      \\\midrule
\multirow{9}{*}{59}   & \multirow{3}{*}{40}     & 256                       & 388.5                          & 22.6                           & 376.1    & 91.0      & 371.1    & 94.9      \\
                      &                         & 512                       & 383.4                          & 25.1                           & 371.1    & 53.5      & 371.1    & 72.3      \\
                      &                         & 2048                     & 371.1                          & 57.0                           & 371.1    & 60.8      & 371.1    & 63.8      \\\cline{2-9} 
                      & \multirow{3}{*}{100}    & 256                       & 388.5                          & 21.3                           & 371.1    & 74.5      & 371.1    & 97.0      \\
                      &                         & 512                       & 388.5                          & 24.1                           & 371.1    & 49.7      & 371.1    & 69.5      \\
                      &                         & 2048                     & 371.1                          & 53.9                           & 371.1    & 61.5      & 371.1    & 66.0      \\\cline{2-9} 
                      & \multirow{3}{*}{400}    & 256                       & 389.5                          & 18.1                           & 371.1    & 73.3      & 371.1    & 95.1      \\
                      &                         & 512                       & 382.0                          & 24.8                           & 371.1    & 54.8      & 371.1    & 71.7      \\
                      &                         & 2048                     & 371.1                          & 55.6                           & 371.1    & 62.4      & 371.1    & 67.6      \\\midrule
\multirow{9}{*}{73}   & \multirow{3}{*}{40}     & 256                       & 338.7                          & 16.4                           & 327.0    & 52.0      & 327.0    & 61.9      \\
                      &                         & 512                       & 342.4                          & 23.9                           & 327.0    & 47.8      & 327.0    & 58.4      \\
                      &                         & 2048                     & 324.7                          & 51.5                           & 324.7    & 59.2      & 324.7    & 65.3      \\\cline{2-9} 
                      & \multirow{3}{*}{100}    & 256                       & 338.7                          & 16.3                           & 328.5    & 48.8      & 328.5    & 58.2      \\
                      &                         & 512                       & 342.4                          & 24.0                           & 328.5    & 36.3      & 328.5    & 43.3      \\
                      &                         & 2048                     & 324.7                          & 55.3                           & 324.7    & 57.8      & 324.7    & 62.7      \\\cline{2-9} 
                      & \multirow{3}{*}{400}    & 256                       & 338.7                          & 18.7                           & 328.5    & 42.8      & 328.5    & 52.2      \\
                      &                         & 512                       & 338.7                          & 25.0                           & 328.5    & 35.9      & 328.5    & 47.3      \\
                      &                         & 2048                     & 324.7                          & 50.2                           & 324.7    & 57.8      & 324.7    & 62.9      \\\bottomrule
\end{tabular}
}
\\[10pt] % Adds some space after the table
\small{
    \emph{search} = width of embedded search, \emph{DD width} = width of relaxed dd \\
    \emph{time} = time to solve
}
\end{table}

\subsection{Second Experiment: Test of Smaller Instances}
\label{sec:second_exp}
\begin{figure}[tb]
  \centering%
  \includegraphics[width=\textwidth]{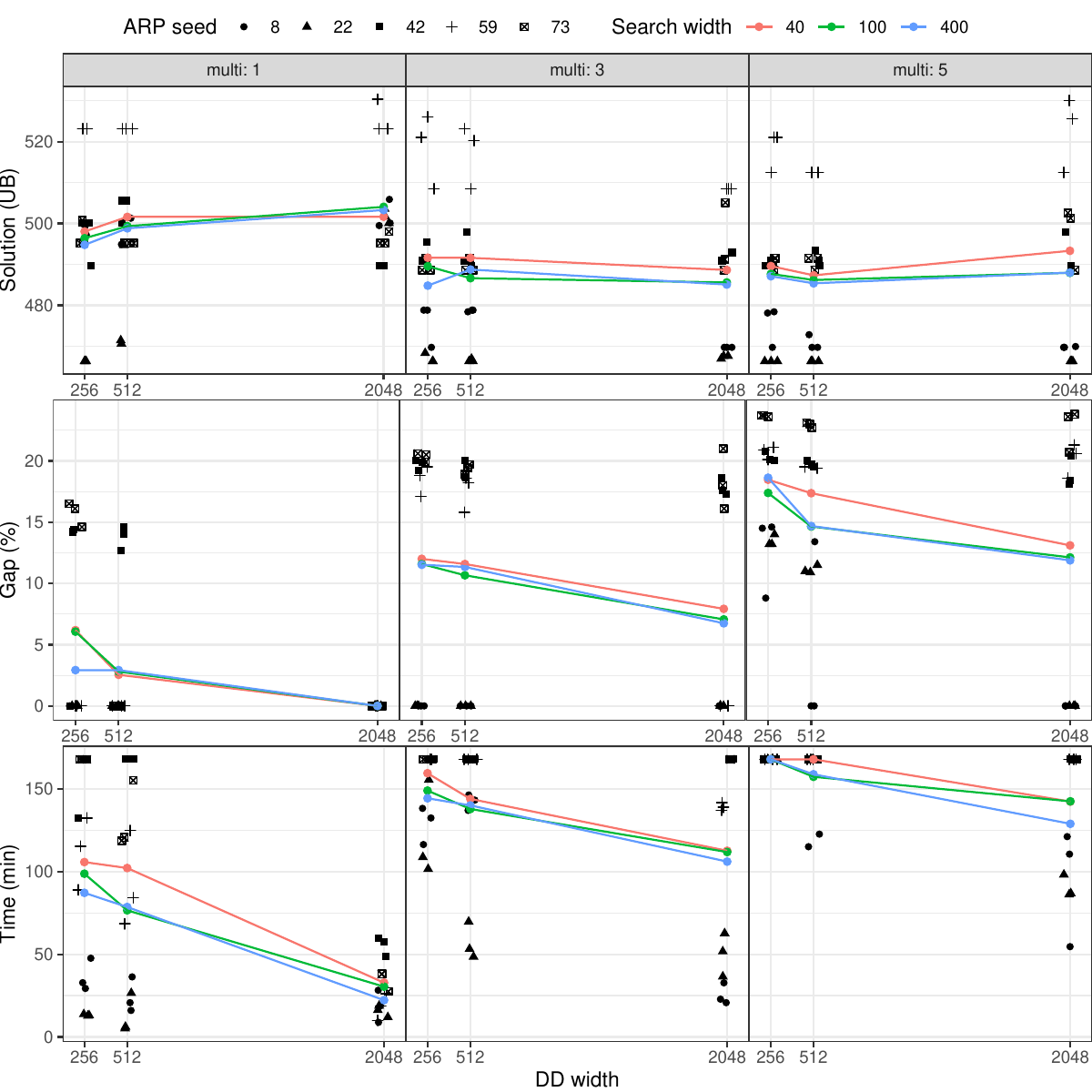}
  \caption{ARP instances with $n=15$ and a 7 day runtime (Table~\ref{tab:test2_15})}
  \label{fig:test2_15}
\end{figure}

\begin{figure}[tb]
  \centering%
  \includegraphics[width=\textwidth]{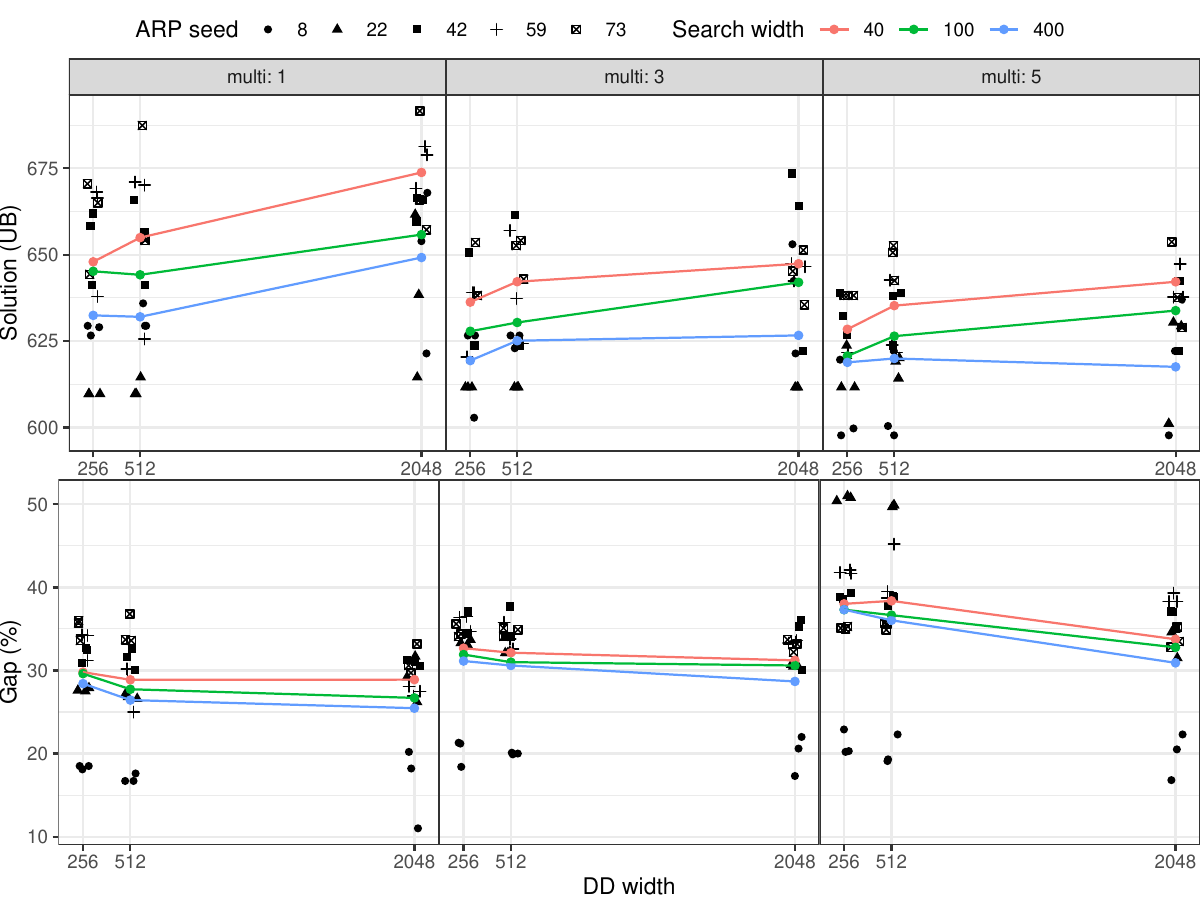}
  \caption{ARP instances with $n=20$ and a 7 day runtime (Table~\ref{tab:test2_20})}
  \label{fig:test2_20}
\end{figure}

Our first round of experiments gave us a rough idea of what settings might be successful. Our second round of experiments makes up the bulk of our tests. As discussed in the previous paragraph we test relaxed DD-widths of 256, 512, and 2048. For the embedded search widths ($\omega_s$ from Algorithm \ref{algo:search}) we tested 40, 100, and 400, where previously we just used 100. In typical DD literature \citep{BerCirHoeHoo2016dd4o}, the size of the search for feasible solutions would be as big or bigger than the size of the relaxed DD, because it is relatively computationally cheap. However, for the ARP, the search for feasible solutions is computationally expensive because of the inner-optimization that requires evaluating  the black-box function  $\mathcal{B}$. This suggested to us that larger embedded search widths are unlikely to be worth it. However, the following results indicate that might not be true. We also tested $\pmulti \in \{1,3,5\}$. Table \ref{tab:test2_10} (Fig.~\ref{fig:test2_10}) show the results of the second round of experiments for $n=10$, Table \ref{tab:test2_15} (Fig.~\ref{fig:test2_15}) for $n=15$, and Table \ref{tab:test2_20} (Fig.~\ref{fig:test2_20}) for $n=20$. In Table \ref{tab:test2_10}, every solution is optimal, so it does not include the \emph{gap} columns.

\begin{table}[tb]
\caption{ARP instances with $n=15$ and a 7 day runtime}\label{tab:test2_15}
\centering
\resizebox{!}{6.75cm}{%
\begin{tabular}{|c|cc|ccccc|ccccc|ccccc|}
\toprule
\multirow{2}{*}{seed} & \multirow{2}{*}{search} & \multirow{2}{*}{DD width} & \multicolumn{5}{c|}{$\pmulti=1$} & \multicolumn{5}{c|}{$\pmulti=3$} & \multicolumn{5}{c|}{$\pmulti=5$}                                                                                        \\
                      &                         &                           & lb                               & ub                               & gap (\%) & time(hr) & queue & lb    & ub    & gap (\%) & time(hr) & queue & lb    & ub    & gap (\%) & time(hr) & queue \\\midrule
\multirow{9}{*}{8}    & \multirow{3}{*}{40}     & 256                       & 499.8                            & 499.8                            & 0.0      & 47.6     & 0     & 478.8 & 478.8 & 0.0      & 138.3    & 0     & 408.1 & 478.1 & 14.6     & -        & 1172  \\
                      &                         & 512                       & 494.9                            & 494.9                            & 0.0      & 36.3     & 0     & 478.4 & 478.4 & 0.0      & 146.4    & 0     & 409.3 & 472.8 & 13.4     & -        & 579   \\
                      &                         & 2048                     & 500.1                            & 500.1                            & 0.0      & 28.2     & 0     & 469.7 & 469.7 & 0.0      & 22.8     & 0     & 469.7 & 469.7 & 0.0      & 121.2    & 0     \\\cline{2-18} 
                      & \multirow{3}{*}{100}    & 256                       & 497.2                            & 497.2                            & 0.0      & 29.3     & 0     & 478.8 & 478.8 & 0.0      & 132.5    & 0     & 428.4 & 469.7 & 8.8      & -        & 311   \\
                      &                         & 512                       & 501.3                            & 501.3                            & 0.0      & 20.7     & 0     & 478.8 & 478.8 & 0.0      & 137.1    & 0     & 469.7 & 469.7 & 0.0      & 115.1    & 0     \\
                      &                         & 2048                     & 505.9                            & 505.9                            & 0.0      & 18.8     & 0     & 469.7 & 469.7 & 0.0      & 32.7     & 0     & 469.9 & 469.9 & 0.0      & 110.6    & 0     \\\cline{2-18} 
                      & \multirow{3}{*}{400}    & 256                       & 499.5                            & 499.5                            & 0.0      & 32.8     & 0     & 469.7 & 469.7 & 0.0      & 116.4    & 0     & 409.2 & 478.4 & 14.5     & -        & 1011  \\
                      &                         & 512                       & 499.5                            & 499.5                            & 0.0      & 16.0     & 0     & 478.8 & 478.8 & 0.0      & 143.3    & 0     & 469.7 & 469.7 & 0.0      & 122.7    & 0     \\
                      &                         & 2048                     & 499.5                            & 499.5                            & 0.0      & 8.8      & 0     & 469.7 & 469.7 & 0.0      & 20.7     & 0     & 469.7 & 469.7 & 0.0      & 54.6     & 0     \\\midrule
\multirow{9}{*}{22}   & \multirow{3}{*}{40}     & 256                       & 466.3                            & 466.3                            & 0.0      & 12.8     & 0     & 466.3 & 466.3 & 0.0      & 155.6    & 0     & 404.6 & 466.3 & 13.2     & -        & 1152  \\
                      &                         & 512                       & 494.7                            & 494.7                            & 0.0      & 26.5     & 0     & 466.9 & 466.9 & 0.0      & 69.7     & 0     & 415.0 & 466.3 & 11.0     & -        & 179   \\
                      &                         & 2048                     & 500.1                            & 500.1                            & 0.0      & 11.9     & 0     & 466.9 & 466.9 & 0.0      & 62.7     & 0     & 466.3 & 466.3 & 0.0      & 86.9     & 0     \\\cline{2-18} 
                      & \multirow{3}{*}{100}    & 256                       & 466.3                            & 466.3                            & 0.0      & 13.0     & 0     & 468.3 & 468.3 & 0.0      & 108.8    & 0     & 404.7 & 466.3 & 13.2     & -        & 1099  \\
                      &                         & 512                       & 471.4                            & 471.4                            & 0.0      & 5.1      & 0     & 466.3 & 466.3 & 0.0      & 48.5     & 0     & 415.6 & 466.3 & 10.9     & -        & 167   \\
                      &                         & 2048                     & 503.6                            & 503.6                            & 0.0      & 18.8     & 0     & 467.6 & 467.6 & 0.0      & 51.7     & 0     & 466.3 & 466.3 & 0.0      & 98.2     & 0     \\\cline{2-18} 
                      & \multirow{3}{*}{400}    & 256                       & 466.3                            & 466.3                            & 0.0      & 13.7     & 0     & 466.3 & 466.3 & 0.0      & 101.6    & 0     & 401.2 & 466.3 & 14.0     & -        & 1058  \\
                      &                         & 512                       & 470.6                            & 470.6                            & 0.0      & 6.1      & 0     & 466.3 & 466.3 & 0.0      & 53.2     & 0     & 412.9 & 466.3 & 11.5     & -        & 206   \\
                      &                         & 2048                     & 501.7                            & 501.7                            & 0.0      & 16.2     & 0     & 467.6 & 467.6 & 0.0      & 36.5     & 0     & 466.3 & 466.3 & 0.0      & 86.2     & 0     \\\midrule
\multirow{9}{*}{42}   & \multirow{3}{*}{40}     & 256                       & 428.0                            & 500.1                            & 14.4     & -        & 1259  & 396.5 & 495.4 & 20.0     & -        & 1672  & 392.6 & 490.9 & 20.0     & -        & 1634  \\
                      &                         & 512                       & 436.7                            & 500.1                            & 12.7     & -        & 316   & 398.2 & 498.0 & 20.0     & -        & 747   & 394.6 & 493.5 & 20.0     & -        & 712   \\
                      &                         & 2048                     & 489.7                            & 489.7                            & 0.0      & 59.7     & 0     & 401.4 & 492.9 & 18.6     & -        & 141   & 396.4 & 497.9 & 20.4     & -        & 147   \\\cline{2-18} 
                      & \multirow{3}{*}{100}    & 256                       & 428.9                            & 500.1                            & 14.2     & -        & 1203  & 396.6 & 490.9 & 19.2     & -        & 1643  & 391.4 & 489.7 & 20.1     & -        & 1537  \\
                      &                         & 512                       & 434.8                            & 505.6                            & 14.0     & -        & 414   & 399.8 & 490.9 & 18.6     & -        & 752   & 395.2 & 490.9 & 19.5     & -        & 709   \\
                      &                         & 2048                     & 489.7                            & 489.7                            & 0.0      & 57.7     & 0     & 406.0 & 490.9 & 17.3     & -        & 105   & 399.6 & 489.7 & 18.4     & -        & 127   \\\cline{2-18} 
                      & \multirow{3}{*}{400}    & 256                       & 489.7                            & 489.7                            & 0.0      & 132.3    & 0     & 392.9 & 490.9 & 19.9     & -        & 1312  & 388.0 & 489.7 & 20.8     & -        & 1232  \\
                      &                         & 512                       & 431.7                            & 505.6                            & 14.6     & -        & 397   & 398.3 & 489.7 & 18.7     & -        & 671   & 393.1 & 489.7 & 19.7     & -        & 646   \\
                      &                         & 2048                     & 489.7                            & 489.7                            & 0.0      & 48.5     & 0     & 404.4 & 490.9 & 17.6     & -        & 109   & 400.8 & 489.7 & 18.1     & -        & 126   \\\midrule
\multirow{9}{*}{59}   & \multirow{3}{*}{40}     & 256                       & 523.2                            & 523.2                            & 0.0      & 132.5    & 0     & 423.4 & 526.1 & 19.5     & -        & 1617  & 412.0 & 521.1 & 20.9     & -        & 1490  \\
                      &                         & 512                       & 523.2                            & 523.2                            & 0.0      & 124.9    & 0     & 427.8 & 523.2 & 18.2     & -        & 721   & 412.9 & 512.5 & 19.4     & -        & 661   \\
                      &                         & 2048                     & 523.2                            & 523.2                            & 0.0      & 26.0     & 0     & 508.5 & 508.5 & 0.0      & 141.7    & 0     & 417.1 & 530.1 & 21.3     & -        & 156   \\\cline{2-18} 
                      & \multirow{3}{*}{100}    & 256                       & 523.2                            & 523.2                            & 0.0      & 115.4    & 0     & 423.1 & 521.1 & 18.8     & -        & 1517  & 411.0 & 521.1 & 21.1     & -        & 1386  \\
                      &                         & 512                       & 523.2                            & 523.2                            & 0.0      & 68.4     & 0     & 428.1 & 508.5 & 15.8     & -        & 655   & 411.7 & 512.5 & 19.7     & -        & 591   \\
                      &                         & 2048                     & 523.2                            & 523.2                            & 0.0      & 18.8     & 0     & 508.5 & 508.5 & 0.0      & 138.9    & 0     & 417.4 & 512.5 & 18.6     & -        & 144   \\\cline{2-18} 
                      & \multirow{3}{*}{400}    & 256                       & 523.2                            & 523.2                            & 0.0      & 89.0     & 0     & 421.5 & 508.5 & 17.1     & -        & 1330  & 409.5 & 512.5 & 20.1     & -        & 1208  \\
                      &                         & 512                       & 523.2                            & 523.2                            & 0.0      & 84.2     & 0     & 423.6 & 520.3 & 18.6     & -        & 543   & 412.6 & 512.5 & 19.5     & -        & 613   \\
                      &                         & 2048                     & 530.4                            & 530.4                            & 0.0      & 10.0     & 0     & 508.5 & 508.5 & 0.0      & 137.2    & 0     & 417.6 & 525.6 & 20.6     & -        & 147   \\\midrule
\multirow{9}{*}{73}   & \multirow{3}{*}{40}     & 256                       & 418.2                            & 500.8                            & 16.5     & -        & 1,253 & 390.9 & 491.5 & 20.5     & -        & 1,653 & 375.7 & 491.5 & 23.6     & -        & 1,557 \\
                      &                         & 512                       & 495.2                            & 495.2                            & 0.0      & 155.2    & 0     & 394.5 & 491.5 & 19.7     & -        & 737   & 378.4 & 491.5 & 23.0     & -        & 689   \\
                      &                         & 2048                     & 495.2                            & 495.2                            & 0.0      & 38.1     & 0     & 399.3 & 505.1 & 21.0     & -        & 152   & 382.9 & 502.6 & 23.8     & -        & 149   \\\cline{2-18} 
                      & \multirow{3}{*}{100}    & 256                       & 415.5                            & 495.2                            & 16.1     & -        & 1,010 & 391.2 & 488.6 & 19.9     & -        & 1,672 & 375.1 & 491.5 & 23.7     & -        & 1,509 \\
                      &                         & 512                       & 495.2                            & 495.2                            & 0.0      & 120.7    & 0     & 396.3 & 488.6 & 18.9     & -        & 741   & 378.0 & 491.5 & 23.1     & -        & 671   \\
                      &                         & 2048                     & 498.0                            & 498.0                            & 0.0      & 37.9     & 0     & 403.0 & 491.2 & 18.0     & -        & 134   & 383.2 & 501.3 & 23.6     & -        & 144   \\\cline{2-18} 
                      & \multirow{3}{*}{400}    & 256                       & 422.7                            & 495.2                            & 14.6     & -        & 1,098 & 387.7 & 488.6 & 20.6     & -        & 1,361 & 372.9 & 488.6 & 23.7     & -        & 1,284 \\
                      &                         & 512                       & 495.2                            & 495.2                            & 0.0      & 118.7    & 0     & 393.9 & 488.6 & 19.4     & -        & 658   & 377.6 & 488.6 & 22.7     & -        & 641   \\
                      &                         & 2048                     & 495.2                            & 495.2                            & 0.0      & 27.5     & 0     & 410.1 & 488.6 & 16.1     & -        & 77    & 387.5 & 488.6 & 20.7     & -        & 135   \\\bottomrule
\end{tabular}
}
\\[10pt] % Adds some space after the table
\small{
    \emph{search} = width of embedded search, \emph{DD width} = width of relaxed dd \\
    \emph{lb} = lower bound, \emph{ub} = upper bound, \emph{gap} = optimality gap,\\
    \emph{time} = time to solve if solved, \emph{queue} = nodes in the processing queue
}
\end{table}

The largest search settings, i.e., an embedded search width of 400 and a relaxed DD-width of 2048, clearly dominate the results, with tighter optimality gaps and more problems solved. For the remainder of the discussion, we will be referring primarily to those results. Recall from Section \ref{sec:limits} that different $\pmulti$ values can have different optimal solutions because the inner optimizer may return a local optimum instead of the global optimum. When $n=10$ (Table \ref{tab:test2_10} and Fig.~\ref{fig:test2_10}) and $\pmulti=1$, all of the problems are solved to optimality in under $1$ hour, and when $\pmulti=5$, they are solved to optimality in under $2$ hours. In some cases the largest search settings alleviates some of the negative effects of $\pmulti=1$ (in Table \ref{tab:test2_10}, seed 22, DD width 512, compare search width 40 to search width 100). This is likely because a larger search provides more opportunities to stumble on improved solutions that will later be trimmed as a result of using a low $\pmulti$. For the the largest search settings, the effect of $\pmulti$ is only to increase the runtime, however settings with smaller DD widths make it clear that $\pmulti$ can have a strong impact on the underlying arc calculations. For example, seed $59$, with search width 40, and relaxed DD-width $256$, has a different exact solution with  $\pmulti=1$ than it does with  $\pmulti=5$; the cost drops from $388.5$ to $371.1$. Recall that because the inner problem in the ARP is nonlinear and nonconvex, it is not possible to prove optimality in general. In lieu of that, $\pmulti$ can be increased arbitrarily to increase confidence in the solutions. If $\pmulti$ is increased, and the algorithm returns the same permutation as it did for a lower $\pmulti$ value, that provides evidence that the optimal solution is stable. The larger the increase in $\pmulti$, the stronger the evidence.

\begin{table}[!ht]
\caption{ARP instances with $n=20$ and a 7 day runtime}\label{tab:test2_20}
\centering
\resizebox{!}{6.75cm}{%
    \begin{tabular}{|c|cc|cccc|cccc|cccc|}
\toprule
\multirow{2}{*}{seed} & \multirow{2}{*}{search} & \multirow{2}{*}{DD width} & \multicolumn{4}{c|}{$\pmulti=1$} & \multicolumn{4}{c|}{$\pmulti=3$} & \multicolumn{4}{c|}{$\pmulti=5$}                                                       \\
                      &                         &                           & lb                               & ub                               & gap (\%) & queue & lb    & ub    & gap (\%) & queue & lb    & ub    & gap (\%) & queue \\\midrule
\multirow{9}{*}{8}    & \multirow{3}{*}{40}     & 256                       & 512.7                            & 629.4                            & 18.5     & 1,012 & 493.5 & 626.6 & 21.2     & 708   & 477.4 & 619.6 & 22.9     & 639   \\
                      &                         & 512                       & 524.0                            & 635.9                            & 17.6     & 557   & 501.4 & 626.6 & 20.0     & 360   & 483.2 & 622.1 & 22.3     & 332   \\
                      &                         & 2048                     & 533.2                            & 667.9                            & 20.2     & 90    & 509.7 & 653.0 & 22.0     & 70    & 495.0 & 637.0 & 22.3     & 70    \\\cline{2-15} 
                      & \multirow{3}{*}{100}    & 256                       & 513.3                            & 626.6                            & 18.1     & 1,101 & 492.9 & 626.6 & 21.3     & 679   & 478.2 & 599.7 & 20.3     & 678   \\
                      &                         & 512                       & 524.2                            & 629.4                            & 16.7     & 598   & 500.4 & 626.6 & 20.1     & 353   & 484.7 & 600.4 & 19.3     & 369   \\
                      &                         & 2048                     & 535.2                            & 653.9                            & 18.2     & 102   & 510.3 & 642.4 & 20.6     & 76    & 494.7 & 622.1 & 20.5     & 75    \\\cline{2-15} 
                      & \multirow{3}{*}{400}    & 256                       & 512.7                            & 629.0                            & 18.5     & 948   & 492.1 & 602.8 & 18.4     & 626   & 476.7 & 597.7 & 20.2     & 581   \\
                      &                         & 512                       & 524.3                            & 629.4                            & 16.7     & 509   & 499.3 & 622.9 & 19.9     & 328   & 483.4 & 597.7 & 19.1     & 344   \\
                      &                         & 2048                     & 552.9                            & 621.4                            & 11.0     & 84    & 514.1 & 621.4 & 17.3     & 94    & 497.5 & 597.7 & 16.8     & 91    \\\midrule
\multirow{9}{*}{22}   & \multirow{3}{*}{40}     & 256                       & 442.3                            & 609.7                            & 27.5     & 710   & 408.8 & 611.6 & 33.1     & 644   & 305.6 & 623.7 & 51.0     & 629   \\
                      &                         & 512                       & 447.5                            & 614.5                            & 27.2     & 380   & 415.6 & 611.6 & 32.1     & 331   & 310.9 & 620.1 & 49.9     & 320   \\
                      &                         & 2048                     & 451.8                            & 661.6                            & 31.7     & 69    & 423.9 & 611.6 & 30.7     & 69    & 410.6 & 630.4 & 34.9     & 66    \\\cline{2-15} 
                      & \multirow{3}{*}{100}    & 256                       & 441.2                            & 609.7                            & 27.6     & 649   & 407.7 & 611.6 & 33.3     & 590   & 303.5 & 611.6 & 50.4     & 578   \\
                      &                         & 512                       & 447.5                            & 609.7                            & 26.6     & 389   & 415.4 & 611.6 & 32.1     & 319   & 310.3 & 619.1 & 49.9     & 314   \\
                      &                         & 2048                     & 453.8                            & 614.5                            & 26.2     & 80    & 425.4 & 611.6 & 30.5     & 69    & 411.5 & 629.3 & 34.6     & 66    \\\cline{2-15} 
                      & \multirow{3}{*}{400}    & 256                       & 439.8                            & 609.7                            & 27.9     & 582   & 405.6 & 611.6 & 33.7     & 515   & 301.3 & 611.6 & 50.8     & 496   \\
                      &                         & 512                       & 446.2                            & 609.7                            & 26.8     & 357   & 414.2 & 611.6 & 32.3     & 297   & 308.8 & 614.2 & 49.7     & 284   \\
                      &                         & 2048                     & 451.8                            & 638.4                            & 29.2     & 69    & 425.4 & 611.6 & 30.5     & 68    & 412.1 & 601.1 & 31.5     & 67    \\\midrule
\multirow{9}{*}{42}   & \multirow{3}{*}{40}     & 256                       & 445.6                            & 661.9                            & 32.7     & 616   & 409.0 & 650.6 & 37.1     & 613   & 387.8 & 638.8 & 39.3     & 598   \\
                      &                         & 512                       & 448.7                            & 665.8                            & 32.6     & 321   & 412.3 & 661.5 & 37.7     & 331   & 391.0 & 638.8 & 38.8     & 331   \\
                      &                         & 2048                     & 457.6                            & 666.5                            & 31.3     & 66    & 430.6 & 673.5 & 36.1     & 65    & 404.3 & 642.3 & 37.1     & 65    \\\cline{2-15} 
                      & \multirow{3}{*}{100}    & 256                       & 444.4                            & 658.2                            & 32.5     & 591   & 409.2 & 623.7 & 34.4     & 597   & 386.9 & 632.2 & 38.8     & 575   \\
                      &                         & 512                       & 448.8                            & 656.5                            & 31.6     & 327   & 411.3 & 623.7 & 34.1     & 301   & 389.1 & 638.0 & 39.0     & 297   \\
                      &                         & 2048                     & 459.3                            & 665.8                            & 31.0     & 67    & 429.8 & 664.2 & 35.3     & 63    & 404.8 & 642.3 & 37.0     & 63    \\\cline{2-15} 
                      & \multirow{3}{*}{400}    & 256                       & 443.1                            & 641.2                            & 30.9     & 548   & 408.6 & 623.7 & 34.5     & 522   & 384.9 & 626.7 & 38.6     & 499   \\
                      &                         & 512                       & 448.4                            & 641.2                            & 30.1     & 290   & 411.0 & 623.7 & 34.1     & 276   & 388.1 & 623.7 & 37.8     & 273   \\
                      &                         & 2048                     & 458.3                            & 659.6                            & 30.5     & 66    & 435.0 & 622.2 & 30.1     & 65    & 403.2 & 622.2 & 35.2     & 63    \\\midrule
\multirow{9}{*}{59}   & \multirow{3}{*}{40}     & 256                       & 439.0                            & 668.1                            & 34.3     & 612   & 406.4 & 639.0 & 36.4     & 594   & 361.9 & 621.7 & 41.8     & 589   \\
                      &                         & 512                       & 468.4                            & 671.0                            & 30.2     & 331   & 422.1 & 657.0 & 35.8     & 316   & 352.4 & 642.6 & 45.2     & 311   \\
                      &                         & 2048                     & 489.6                            & 681.3                            & 28.1     & 68    & 429.9 & 647.4 & 33.6     & 70    & 392.7 & 647.3 & 39.3     & 67    \\\cline{2-15} 
                      & \multirow{3}{*}{100}    & 256                       & 438.5                            & 666.4                            & 34.2     & 592   & 405.8 & 639.0 & 36.5     & 568   & 359.7 & 621.6 & 42.1     & 560   \\
                      &                         & 512                       & 468.2                            & 670.1                            & 30.1     & 325   & 421.9 & 637.3 & 33.8     & 308   & 377.4 & 623.9 & 39.5     & 301   \\
                      &                         & 2048                     & 492.6                            & 678.9                            & 27.5     & 71    & 430.7 & 646.5 & 33.4     & 68    & 393.3 & 637.7 & 38.3     & 69    \\\cline{2-15} 
                      & \multirow{3}{*}{400}    & 256                       & 438.9                            & 637.9                            & 31.2     & 595   & 404.8 & 620.4 & 34.7     & 514   & 361.7 & 619.9 & 41.7     & 469   \\
                      &                         & 512                       & 469.2                            & 625.6                            & 25.0     & 351   & 420.9 & 624.3 & 32.6     & 293   & 381.1 & 621.6 & 38.7     & 283   \\
                      &                         & 2048                     & 489.3                            & 669.2                            & 26.9     & 72    & 428.4 & 642.4 & 33.3     & 67    & 393.4 & 637.7 & 38.3     & 66    \\\midrule
\multirow{9}{*}{73}   & \multirow{3}{*}{40}     & 256                       & 429.0                            & 670.5                            & 36.0     & 630   & 420.8 & 653.5 & 35.6     & 614   & 414.0 & 638.2 & 35.1     & 559   \\
                      &                         & 512                       & 434.1                            & 687.4                            & 36.8     & 312   & 424.6 & 654.1 & 35.1     & 314   & 419.5 & 652.6 & 35.7     & 313   \\
                      &                         & 2048                     & 462.0                            & 691.6                            & 33.2     & 74    & 431.6 & 651.3 & 33.7     & 68    & 423.9 & 653.7 & 35.2     & 66    \\\cline{2-15} 
                      & \multirow{3}{*}{100}    & 256                       & 427.8                            & 665.0                            & 35.7     & 562   & 420.3 & 638.2 & 34.1     & 571   & 415.0 & 638.2 & 35.0     & 558   \\
                      &                         & 512                       & 434.6                            & 655.1                            & 33.7     & 321   & 424.6 & 652.6 & 34.9     & 316   & 419.4 & 650.6 & 35.6     & 309   \\
                      &                         & 2048                     & 462.2                            & 665.9                            & 30.6     & 69    & 430.7 & 645.2 & 33.2     & 66    & 423.9 & 637.6 & 33.5     & 68    \\\cline{2-15} 
                      & \multirow{3}{*}{400}    & 256                       & 427.5                            & 644.3                            & 33.6     & 543   & 418.9 & 638.2 & 34.4     & 485   & 412.9 & 638.2 & 35.3     & 507   \\
                      &                         & 512                       & 434.1                            & 654.1                            & 33.6     & 300   & 423.8 & 642.9 & 34.1     & 291   & 418.2 & 642.5 & 34.9     & 288   \\
                      &                         & 2048                     & 462.2                            & 657.2                            & 29.7     & 70    & 431.1 & 635.5 & 32.2     & 63    & 422.7 & 628.9 & 32.8     & 65    \\\bottomrule
\end{tabular}
}
\\[10pt] % Adds some space after the table
\small{
    \emph{search} = width of embedded search, \emph{DD width} = width of relaxed dd\\
    \emph{lb} = lower bound, \emph{ub} = upper bound, \emph{gap} = optimality gap,\\
    \emph{queue} = nodes in the processing queue
}
\end{table}

We now focus on the results obtained with  an embedded search width of 400 and a relaxed DD-width of 2048 on instances of size $n \in \{15,20\}$ (Tables \ref{tab:test2_15} and \ref{tab:test2_20}) with a maximum runtime of 7 days. With $\pmulti=1$, three of the problems were solved in under a day, with the other two taking around two days. However, many of the solutions are improved with a higher $\pmulti$ (see the top row of Figures \ref{fig:test2_15} and \ref{fig:test2_20}). For $n=15$, the setting $\pmulti=1$ is unlikely to yield optimal solutions. With $\pmulti=5$, only seeds $8$ and $22$ were solved exactly, but they yielded much better solutions. In seed $22$ for example, the best solution drops from $501.7$ to $466.3$. With $n=20$, none of the problems were solved to optimality within the 7 day limit.

The only available comparison we can make is with \cite{LopChiGil2022evo}. However, because they assume that the inner problem is a total black box, and we take advantage of the structure of the black box, any time/computation comparison is meaningless. So we will merely be comparing our solutions to theirs, to show that our method yields high quality solutions. They only used seeds 42 and 73, so for the other problems we can only report our solutions. These results are shown in Table \ref{tab:test2_c}. Notably the solution reported by Peel-and-Bound's preferred setting ($\omega_s = 400$, $\omega = 2048$, \pmulti$=5$) is the same as the best previously found for $n=10$, and better than those previously reported for $n=15$ and $n=20$. In some cases, the gap is quite large, for example with $n=20$, seed$=42$, our solution is $12.0\%$ better. We also report the best observed solution from any setting, and for 13 of the 15 instances the preferred setting found the same solution as the best solution found by any setting. All of theses results use $\mathcal{B}$ to evaluate the cost of a permutation, and thus they are directly comparable. In other words, the values listed are not dependent on the optimality of the inner problem.

\begin{table}[!ht]
\caption{Best solutions for $n \in \{10,15,20\}$}
\centering
\resizebox{!}{2.3cm}{%
\begin{tabular}{|c|c|cc|cc|cc|}
\hline
\multirow{2}{*}{n} & \multirow{2}{*}{seed} & \multicolumn{2}{c|}{\cite{LopChiGil2022evo}$^*$} & \multicolumn{2}{c|}{Peel-and-Bound: Preferred Settings$^{**}$} & \multicolumn{2}{c|}{Peel-and-Bound: Best Found$^{***}$} \\
 &  & \multicolumn{1}{l}{Average Value} & Best Value & Value & Optimality Gap ($\%$) & Value & Optimality Gap ($\%$) \\ \hline
\multirow{5}{*}{10} & 8 & \multicolumn{1}{c|}{-} & - & 360.4 & 0.0 & 357.5 & $0.0$ \\
 & 22 & \multicolumn{1}{c|}{-} & - & 364.5 & 0.0 & 364.5 & $0.0$ \\
 & 42 & \multicolumn{1}{c|}{374.9} & 346.7 & 346.7 & 0.0 & 346.7 & $0.0$ \\
 & 59 & \multicolumn{1}{c|}{-} & - & 371.1 & 0.0 & 371.1 & $0.0$ \\
 & 73 & \multicolumn{1}{c|}{355.9} & 324.7 & 324.7 & 0.0 & 324.7 & $0.0$ \\ \hline
\multirow{5}{*}{15} & 8 & \multicolumn{1}{c|}{-} & - & 469.7 & 0.0 & 469.7 & $0.0$ \\
 & 22 & \multicolumn{1}{c|}{-} & - & 466.3 & 0.0 & 466.3 & $0.0$ \\
 & 42 & \multicolumn{1}{c|}{497.2} & 490.9 & 489.7 & 18.1 & 489.7 & $0.0$ \\
 & 59 & \multicolumn{1}{c|}{-} & - & 525.6 & 20.6 & \textbf{508.5} & $0.0$ \\
 & 73 & \multicolumn{1}{c|}{525.6} & 519.9 & 488.6 & 20.7 & 488.6 & 16.1 \\ \hline
\multirow{5}{*}{20} & 8 & \multicolumn{1}{c|}{-} & - & 597.7 & 16.8 & 597.7 & 20.3 \\
 & 22 & \multicolumn{1}{c|}{-} & - & 601.1 & 31.5 & 601.1 & 31.5 \\
 & 42 & \multicolumn{1}{c|}{737.0} & 707.2 & 622.2 & 35.2 & 622.2 & 30.1 \\
 & 59 & \multicolumn{1}{c|}{-} & - & 637.7 & 38.3 & \textbf{619.9} & 41.7 \\
 & 73 & \multicolumn{1}{c|}{661.8} & 652.5 & 628.9 & 32.8 & 628.9 & 32.8 \\ \hline
\end{tabular}
}
\\[10pt] % Adds some space after the table
\small{
    $^*$ The results from  \cite{LopChiGil2022evo} in the \emph{Average Found} column are the average solution cost over several runs of the same stochastic algorithm. \emph{Best Found} reports the best solution found in any run of the algorithm.\\
    $^{**}$ The preferred settings for Peel-and-Bound are: $\omega_s = 400$, $\omega = 2048$, \pmulti$= 5$. \\
    $^{***}$ The bold values are the only instances where the best found solution by any setting is better than the solution found using the preferred setting.
}
\label{tab:test2_c}

\end{table}

\subsection{Final Experiment: Test of Larger Instances}

In our final round of experiments we run the larger instances with $n \in \{25, 30\}$. Due to limited resource availability we had to use a slightly less powerful computer equipped with an Intel E5-2650 v4 Broadwell 2.2GHz CPU with 64Gb RAM. We use a relaxed DD-width of 2048, and we test larger embedded search sizes than before: 400, 1024, and 2048. We limit the runtime to 3 days, again due to resource availability, but this has the added benefit of showing that the algorithm works well within a smaller time-frame than the 7 days used in Section \ref{sec:second_exp}.

The results, shown in Tables \ref{tab:test3_25} and \ref{tab:test3_30}, suggest that increasing the width of the embedded search may be worth the extra computational cost of constructing the additional solutions, but the results are not enough to be conclusive. $\omega_s = 2048$ finds the best solutions more often than the other settings, so as before we will choose the largest search setting to be the preferred setting. More experiments would be needed to make a strong recommendation between the search settings, but it is clear that all of the settings used in this section are effective.

\begin{table}[!ht]
\caption{ARP instances with $n=25$, a relaxed DD-width of 2048, and a 3 day runtime}
\centering
\resizebox{!}{2.5cm}{%
\begin{tabular}{|c|c|cccc|cccc|cccc|}
\midrule
\multirow{2}{*}{Seed} & \multirow{2}{*}{Search} & \multicolumn{4}{c|}{$\pmulti=1$} & \multicolumn{4}{c|}{$\pmulti=3$} & \multicolumn{4}{c|}{$\pmulti=5$}                                                       \\
                      &                         & lb                               & ub                               & gap (\%) & queue & lb    & ub    & gap (\%) & queue & lb    & ub    & gap (\%) & queue \\\midrule
\multirow{3}{*}{8}    & 400                     & 538.1                            & 781.5                            & 31.2     & 4     & 510.3 & 778.7 & 34.5     & 10    & 503.2 & 760.6 & 33.8     & 9     \\
                      & 1024                   & 538.1                            & 778.9                            & 30.9     & 2     & 510.1 & 780.1 & 34.6     & 9     & 503.2 & 760.6 & 33.8     & 9     \\
                      & 2048                   & 538.1                            & 778.9                            & 30.9     & 2     & 510.3 & 774.0 & 34.1     & 9     & 499.7 & 763.3 & 34.5     & 7     \\\midrule
\multirow{3}{*}{22}   & 400                     & 506.8                            & 842.9                            & 39.9     & 10    & 475.4 & 774.4 & 38.6     & 10    & 463.8 & 791.2 & 41.4     & 11    \\
                      & 1024                   & 506.8                            & 787.0                            & 35.6     & 11    & 475.4 & 774.0 & 38.6     & 10    & 463.7 & 785.3 & 41.0     & 7     \\
                      & 2048                   & 506.8                            & 784.6                            & 35.4     & 10    & 475.4 & 774.0 & 38.6     & 9     & 463.7 & 784.6 & 40.9     & 7     \\\midrule
\multirow{3}{*}{42}   & 400                     & 513.9                            & 783.5                            & 34.4     & 2     & 497.5 & 761.1 & 34.6     & 9     & 467.2 & 733.3 & 36.3     & 9     \\
                      & 1024                   & 513.9                            & 776.3                            & 33.8     & 2     & 497.5 & 755.0 & 34.1     & 9     & 467.2 & 751.3 & 37.8     & 9     \\
                      & 2048                   & 513.9                            & 778.7                            & 34.0     & 2     & 497.5 & 755.0 & 34.1     & 9     & 467.2 & 751.3 & 37.8     & 9     \\\midrule
\multirow{3}{*}{59}   & 400                     & 508.6                            & 794.6                            & 36.0     & 6     & 486.0 & 747.6 & 35.0     & 10    & 466.8 & 725.2 & 35.6     & 9     \\
                      & 1024                   & 508.8                            & 795.6                            & 36.1     & 11    & 486.0 & 751.5 & 35.3     & 10    & 466.1 & 724.5 & 35.7     & 7     \\
                      & 2048                   & 508.6                            & 781.3                            & 34.9     & 6     & 481.3 & 738.2 & 34.8     & 7     & 466.8 & 718.2 & 35.0     & 9     \\\midrule
\multirow{3}{*}{73}   & 400                     & 494.7                            & 797.3                            & 38.0     & 9     & 461.2 & 793.8 & 41.9     & 11    & 447.4 & 743.4 & 39.8     & 7     \\
                      & 1024                   & 493.6                            & 790.3                            & 37.5     & 7     & 459.4 & 784.9 & 41.5     & 8     & 447.4 & 743.4 & 39.8     & 7     \\
                      & 2048                   & 493.6                            & 782.0                            & 36.9     & 7     & 459.1 & 768.0 & 40.2     & 7     & 447.4 & 743.4 & 39.8     & 7     \\\midrule
\end{tabular}
}
\\[10pt] % Adds some space after the table
\small{
    \emph{search} = width of embedded search,
    \emph{lb} = lower bound, \emph{ub} = upper bound, 
    \\
    \emph{gap} = optimality gap,
    \emph{queue} = nodes in the processing queue
}
\label{tab:test3_25}
\end{table}

\begin{table}[!ht]
\caption{ARP instances with $n=30$, a relaxed DD-width of 2048, and a 3 day runtime}
\centering
\resizebox{!}{2.5cm}{%
\begin{tabular}{|c|c|cccc|cccc|cccc|}
\toprule
\multirow{2}{*}{Seed} & \multirow{2}{*}{Search} & \multicolumn{4}{c|}{$\pmulti=1$} & \multicolumn{4}{c|}{$\pmulti=3$} & \multicolumn{4}{c|}{$\pmulti=5$}                                                       \\
                      &                         & lb                               & ub                               & gap (\%) & queue & lb    & ub    & gap (\%) & queue & lb    & ub    & gap (\%) & queue \\\midrule
\multirow{3}{*}{8}    & 400                     & 523.5                            & 955.5                            & 45.2     & 4     & 486.9 & 906.9 & 46.3     & 7     & 471.8 & 904.2 & 47.8     & 7     \\
                      & 1024                   & 523.5                            & 935.7                            & 44.1     & 4     & 486.9 & 894.8 & 45.6     & 7     & 464.1 & 898.7 & 48.4     & 6     \\
                      & 2048                   & 523.5                            & 935.7                            & 44.1     & 4     & 486.9 & 892.2 & 45.4     & 7     & 464.1 & 898.7 & 48.4     & 6     \\\midrule
\multirow{3}{*}{22}   & 400                     & 608.8                            & 901.8                            & 32.5     & 5     & 597.5 & 896.2 & 33.3     & 8     & 503.7 & 908.5 & 44.6     & 7     \\
                      & 1024                   & 608.8                            & 896.7                            & 32.1     & 5     & 593.5 & 889.6 & 33.3     & 7     & 503.7 & 905.6 & 44.4     & 7     \\
                      & 2048                   & 608.8                            & 897.9                            & 32.2     & 5     & 597.5 & 900.6 & 33.6     & 6     & 503.5 & 892.8 & 43.6     & 5     \\\midrule
\multirow{3}{*}{42}   & 400                     & 515.9                            & 906.6                            & 43.1     & 9     & 486.7 & 880.1 & 44.7     & 7     & 465.7 & 858.6 & 45.8     & 6     \\
                      & 1024                   & 513.7                            & 911.8                            & 43.7     & 8     & 486.7 & 873.1 & 44.3     & 7     & 465.7 & 854.8 & 45.5     & 6     \\
                      & 2048                   & 512.9                            & 861.1                            & 40.4     & 7     & 486.7 & 871.3 & 44.1     & 7     & 465.7 & 835.1 & 44.2     & 6     \\\midrule
\multirow{3}{*}{59}   & 400                     & 552.0                            & 923.7                            & 40.2     & 6     & 537.4 & 906.2 & 40.7     & 7     & 449.1 & 872.7 & 48.5     & 3     \\
                      & 1024                   & 552.0                            & 910.6                            & 39.4     & 6     & 535.7 & 858.1 & 37.6     & 6     & 449.1 & 864.5 & 48.1     & 3     \\
                      & 2048                   & 552.0                            & 894.8                            & 38.3     & 6     & 534.6 & 868.3 & 38.4     & 2     & 449.1 & 864.1 & 48.0     & 3     \\\midrule
\multirow{3}{*}{73}   & 400                     & 504.1                            & 896.6                            & 43.8     & 5     & 495.1 & 903.9 & 45.2     & 4     & 480.2 & 904.3 & 46.9     & 5     \\
                      & 1024                   & 504.1                            & 975.2                            & 48.3     & 5     & 495.1 & 932.8 & 46.9     & 3     & 480.2 & 906.3 & 47.0     & 5     \\
                      & 2048                   & 504.1                            & 945.1                            & 46.7     & 5     & 495.1 & 908.5 & 45.5     & 2     & 478.5 & 882.9 & 45.8     & 2     \\\bottomrule
\end{tabular}
}
\\[10pt] % Adds some space after the table
\small{
    \emph{search} = width of embedded search,
    \emph{lb} = lower bound, \emph{ub} = upper bound, 
    \\
    \emph{gap} = optimality gap,
    \emph{queue} = nodes in the processing queue
}
\label{tab:test3_30}
\end{table}

We again compare our results with \cite{LopChiGil2022evo} where possible, and record the best solution found using our preferred setting as well as the best solutions found by any setting. This is shown in Table \ref{tab:test3_c}. As before, the solutions from our preferred setting are better than the best solutions found by \cite{LopChiGil2022evo} for all of the problems. The preferred solution only matches the best found solution for 5 of the 10 problems. However, the gap between the two solutions is quite small for the other 5 instances. The largest difference, which occurs with $n=25$ and seed$=42$, is only $2.4\%$.

\FloatBarrier

\begin{table}[!ht]
\caption{Best solutions for $n \in \{25,30\}$}
\centering
\resizebox{!}{2.3cm}{%
\begin{tabular}{|c|c|cc|cc|cc|}
\hline
\multirow{2}{*}{n} & \multirow{2}{*}{seed} & \multicolumn{2}{c|}{\cite{LopChiGil2022evo}$^*$} & \multicolumn{2}{c|}{Peel-and-Bound: Preferred Settings$^{**}$} & \multicolumn{2}{c|}{Peel-and-Bound: Best Found$^{***}$} \\
 &  & Average Value & Best Value & Value & Optimality Gap ($\%$) & Value & Optimality Gap ($\%$) \\ \hline
\multirow{5}{*}{25} & 8 & - & - & 763.3 & 34.5 & \textbf{760.6} & 33.8 \\
 & 22 & - & - & 784.6 & 40.9 & \textbf{774.0} & 38.6 \\
 & 42 & 881.5 & 865.7 & 751.3 & 37.8 & \textbf{733.3} & 36.3 \\
 & 59 & - & - & 718.2 & 35 & 718.2 & 35.0 \\
 & 73 & 873.6 & 863.7 & 743.4 & 39.8 & 743.4 & 39.8 \\ \hline
\multirow{5}{*}{30} & 8 & - & - & 898.7 & 48.4 & \textbf{892.2} & 45.4 \\
 & 22 & - & - & 892.8 & 43.6 & 892.8 & 43.6 \\
 & 42 & 1084.6 & 1065.2 & 835.1 & 44.2 & 835.1 & 44.2 \\
 & 59 & - & - & 864.1 & 48.0 & \textbf{858.1} & 37.6 \\
 & 73 & 967.7 & 952.1 & 882.9 & 45.8 & 882.9 & 45.8 \\ \hline
\end{tabular}
}
\\[10pt] % Adds some space after the table
\small{
    $^*$ The results from  \cite{LopChiGil2022evo} in the \emph{Average Found} column are the average solution cost over several runs of the same stochastic algorithm. \emph{Best Found} reports the best solution found in any run of the algorithm.\\
    $^{**}$ The preferred settings for Peel-and-Bound are: $\omega_s = 2048$, $\omega = 2048$, \pmulti$= 5$. \\
    $^{***}$ The bold values are the instances where the best found solution by any setting is better than the solution found using the preferred setting.
}
\label{tab:test3_c}
\end{table}

\section{Future Work: Opportunities with Parallel Computing}
\label{sec:future}
The process of solving a DD using Peel-and-Bound is extremely parallelizable \cite{BnB,Perez_Regin_2018,RudCapRou2023improved_pnb}. Each peeled DD can be considered a discrete problem to be solved, allowing the problem to be easily divided among any number of available processors, without requiring the processors to communicate (with the exception of requesting new tasks). Consider the number of diagrams left in the queues from the results tables. Large values indicate that: (1) a significant amount of work remains to be done before closing those instances, and (2) that work is extremely parallelizable because each DD in the queue represents a sub-problem to be solved that is totally independent from the rest. 

We ran all of our experiments on a single thread, but as long as there are unused threads available, each of those DDs could be assigned to a different thread without the need for them to communicate. Thus, an implementation that included parallel processing could handle problems at a much larger scale if many CPUs were available. Additionally, setting $\pmulti > 1$ is parallelizable, as it represents $\pmulti$ independent restarts of SQSLP at disjoint areas of the search space.

\section[]{Conclusions and Looking Forward}
\label{sec:conclusion}

In this paper, we study  outer-inner optimization problems, where the outer problem is combinatorial and the inner problem is numerical and black-box. We introduce the first optimization framework designed to find exact solutions for such problems under mild assumptions about the inner problem. Global trajectory optimization is one real-world example of such problems, and we use the Asteroid Routing Problem (ARP) as a case study. Although the inner optimizer in the ARP returns a local optimum, we show how to control the likelihood that this local optimum is also global. From the perspective of the outer problem, our proposed method successfully solves instances of up to 15 celestial bodies, and finds high-quality solutions for larger problems. Notably, we find new best-known solutions for several instances, many of them likely to be optimal. Additionally, we have made our implementation, data, and a detailed guide on how to use the solver, available in a public repository.

In the domain of global trajectory optimization, our methods represent a pioneering approach not only for finding optimal solutions but also for discovering high-quality feasible solutions.  This work opens the door to future research that could adapt our solver to tackle other complex challenges, such as those involving multiple gravity assists. More broadly, we provide a robust framework for addressing sequencing problems where the cost function is computationally demanding. The proposed approach is highly scalable and future work should explore how to set its parameters according to the number of CPUs available and problem size. The principles underlying our approach are versatile, promising applicability to a diverse range of problems in the future.

\ACKNOWLEDGMENT{
\label{sec:acknowledge}
We received advice on this project from the Advanced Concepts Team at the European Space Agency. We would like to specifically thank Dario Izzo and Emmanuel Blazquez for their input on this project.
}

\FloatBarrier
\bibliographystyle{informs2014}
\bibliography{bib/abbrev,bib/journals,bib/authors,bib/articles,bib/biblio,bib/crossref,arp}

\providecommand{\MaxMinAntSystem}{{$\cal MAX$--$\cal MIN$} {Ant} {System}} \providecommand{\rpackage}[1]{{#1}} \providecommand{\softwarepackage}[1]{{#1}} \providecommand{\proglang}[1]{{#1}}
\begin{thebibliography}{34}
\providecommand{\natexlab}[1]{#1}
\providecommand{\url}[1]{\texttt{#1}}
\providecommand{\urlprefix}{URL }

\bibitem[{Abdelkhalik \protect\BIBand{} Gad(2012)}]{AbdGad2012dynamic}
Abdelkhalik O, Gad A (2012) Dynamic-size multiple populations genetic algorithm for multigravity-assist trajectory optimization. \emph{Journal of Guidance, Control, and Dynamics} 35(2):520--529, \urlprefix\url{http://dx.doi.org/10.2514/1.54330}.

\bibitem[{Andersson et~al.(2015)Andersson, Fagerholt, \protect\BIBand{} Hobbesland}]{AndFagHob2015maritime}
Andersson H, Fagerholt K, Hobbesland K (2015) Integrated maritime fleet deployment and speed optimization: Case study from {RoRo} shipping. \emph{Computers \& Operations Research} 55:233--240, \urlprefix\url{http://dx.doi.org/10.1016/j.cor.2014.03.017}.

\bibitem[{Bergman et~al.(2014)Bergman, Cire, Sabharwal, Samulowitz, Saraswat, \protect\BIBand{} van Hoeve}]{BnB}
Bergman D, Cire AA, Sabharwal A, Samulowitz H, Saraswat V, van Hoeve WJ (2014) Parallel combinatorial optimization with decision diagrams. \emph{Proceedings of the International Conference on AI and OR Techniques in Constraint Programming for Combinatorial Optimization Problems}, 351--367.

\bibitem[{Bergman et~al.(2016)Bergman, Cire, van Hoeve, \protect\BIBand{} Hooker}]{BerCirHoeHoo2016dd4o}
Bergman D, Cire AA, van Hoeve WJ, Hooker J (2016) \emph{Decision Diagrams for Optimization} (Cham, Switzerland: Springer), ISBN 978-3-319-42849-9, \urlprefix\url{http://dx.doi.org/10.1007/978-3-319-42849-9}.

\bibitem[{Ceriotti \protect\BIBand{} Vasile(2010)}]{CerVas2010mga}
Ceriotti M, Vasile M (2010) Automated multigravity assist trajectory planning with a modified ant colony algorithm. \emph{Journal of Aerospace Computing, Information, and Communication} 7(9):261--293, \urlprefix\url{http://dx.doi.org/10.2514/1.48448}.

\bibitem[{Chicano et~al.(2023)Chicano, Derbel, \protect\BIBand{} Verel}]{ChiDerVer2023fourier}
Chicano F, Derbel B, Verel S (2023) Fourier transform-based surrogates for permutation problems. Silva S, Paquete L, eds., \emph{Proceedings of the Genetic and Evolutionary Computation Conference, GECCO 2023}, 275--283 (New York, NY: ACM Press), \urlprefix\url{http://dx.doi.org/10.1145/3583131}.

\bibitem[{Cire \protect\BIBand{} van Hoeve(2013)}]{CirHoe2013mdd}
Cire AA, van Hoeve WJ (2013) Multivalued decision diagrams for sequencing problems. \emph{Operations Research} 61(6):1259--1462, \urlprefix\url{http://dx.doi.org/10.1287/opre.2013.1221}.

\bibitem[{Copp{\'e} et~al.(2024)Copp{\'e}, Gillard, \protect\BIBand{} Schaus}]{CopGilSch2024decision}
Copp{\'e} V, Gillard X, Schaus P (2024) Decision diagram-based branch-and-bound with caching for dominance and suboptimality detection. \emph{INFORMS Journal on Computing} \urlprefix\url{http://dx.doi.org/10.1287/ijoc.2022.0340}.

\bibitem[{Ehmke et~al.(2016)Ehmke, Campbell, \protect\BIBand{} Thomas}]{EhmCamTho2016tdvrp}
Ehmke JF, Campbell AM, Thomas BW (2016) Vehicle routing to minimize time-dependent emissions in urban areas. \emph{European Journal of Operational Research} 251(2):478--494, \urlprefix\url{http://dx.doi.org/10.1016/j.ejor.2015.11.034}.

\bibitem[{Gendreau et~al.(2015)Gendreau, Ghiani, \protect\BIBand{} Guerriero}]{GenGhiGue2015tdtsp}
Gendreau M, Ghiani G, Guerriero E (2015) Time-dependent routing problems: A review. \emph{Computers \& Operations Research} 64:189--197, \urlprefix\url{http://dx.doi.org/10.1016/j.cor.2015.06.001}.

\bibitem[{Gillard(2022)}]{Gillard2022phd}
Gillard X (2022) \emph{Discrete Optimization with Decision Diagrams: Design of a Generic Solver, Improved Bounding Techniques, and Discovery of Good Feasible Solutions with Large Neighborhood Search}. Ph.D. thesis, Universit{\'e} Catholique de Louvain.

\bibitem[{Gillard et~al.(2021)Gillard, Copp{\'e}, Schaus, \protect\BIBand{} Cire}]{GilCopSch2021bbmdd}
Gillard X, Copp{\'e} V, Schaus P, Cire AA (2021) Improving the filtering of {Branch-And-Bound} {MDD} solver. Stuckey PJ, ed., \emph{Integration of Constraint Programming, Artificial Intelligence, and Operations Research, CPAIOR 2021}, volume 12735 of \emph{Lecture Notes in Computer Science}, 231--247 (Cham, Switzerland: Springer), \urlprefix\url{http://dx.doi.org/10.1007/978-3-030-78230-6_15}.

\bibitem[{Hansen \protect\BIBand{} Ostermeier(2001)}]{HanOst2001ec}
Hansen N, Ostermeier A (2001) Completely derandomized self-adaptation in evolution strategies. \emph{Evolutionary Computation} 9(2):159--195, \urlprefix\url{http://dx.doi.org/10.1162/106365601750190398}.

\bibitem[{Hennes \protect\BIBand{} Izzo(2015)}]{HenIzz2015interplanetary}
Hennes D, Izzo D (2015) Interplanetary trajectory planning with {Monte} {Carlo} tree search. Yang Q, Wooldridge M, eds., \emph{Proceedings of the 24th International Joint Conference on Artificial Intelligence (IJCAI-15)}, 769--775 (IJCAI/AAAI Press, Menlo Park, CA).

\bibitem[{Hennes et~al.(2016)Hennes, Izzo, \protect\BIBand{} Landau}]{HenIzzLan2016fast}
Hennes D, Izzo D, Landau D (2016) Fast approximators for optimal low-thrust hops between main belt asteroids. Chen X, Stafylopatis A, eds., \emph{Computational Intelligence (SSCI), 2016 IEEE Symposium Series on}, 1--7, \urlprefix\url{http://dx.doi.org/10.1109/SSCI.2016.7850107}.

\bibitem[{Irurozki \protect\BIBand{} L{\'o}pez-Ib{\'a}{\~n}ez(2021)}]{IruLop2021gecco}
Irurozki E, L{\'o}pez-Ib{\'a}{\~n}ez M (2021) Unbalanced mallows models for optimizing expensive black-box permutation problems. Chicano F, Krawiec K, eds., \emph{Proceedings of the Genetic and Evolutionary Computation Conference, GECCO 2021}, 225--233 (New York, NY: ACM Press), \urlprefix\url{http://dx.doi.org/10.1145/3449639.3459366}.

\bibitem[{Izzo(2015)}]{Izzo2015lambert}
Izzo D (2015) Revisiting {Lambert}'s problem. \emph{Celestial Mechanics and Dynamical Astronomy} 121:1--15.

\bibitem[{Izzo et~al.(2007)Izzo, Becerra, Myatt, Nasuto, \protect\BIBand{} Bishop}]{IzzBecMyaNas2007search}
Izzo D, Becerra VM, Myatt DR, Nasuto SJ, Bishop JM (2007) Search space pruning and global optimisation of multiple gravity assist spacecraft trajectories. \emph{Journal of Global Optimization} 38:283--296, \urlprefix\url{http://dx.doi.org/10.1007/s10898-006-9106-0}.

\bibitem[{Izzo et~al.(2015)Izzo, Getzner, Hennes, \protect\BIBand{} Sim{\~o}es}]{IzzGetHenSim2015evolving}
Izzo D, Getzner I, Hennes D, Sim{\~o}es LF (2015) Evolving solutions to {TSP} variants for active space debris removal. Silva S, Esparcia{-}Alc{\'{a}}zar AI, eds., \emph{Proceedings of the Genetic and Evolutionary Computation Conference, GECCO 2015}, 1207--1214 (New York, NY: ACM Press).

\bibitem[{Izzo et~al.(2013)Izzo, Sim{\~o}es, M{\"a}rtens, de~Croon, Heritier, \protect\BIBand{} Yam}]{IzzSimMar2013tour}
Izzo D, Sim{\~o}es LF, M{\"a}rtens M, de~Croon GC, Heritier A, Yam CH (2013) Search for a grand tour of the {Jupiter} {Galilean} moons. Blum C, Alba E, eds., \emph{Proceedings of the Genetic and Evolutionary Computation Conference, GECCO 2013}, 1301--1308 (New York, NY: ACM Press), ISBN 978-1-4503-1963-8, \urlprefix\url{http://dx.doi.org/10.1145/2463372.2463524}.

\bibitem[{Kraft(1988)}]{Kraft1988slsqp}
Kraft D (1988) A software package for sequential quadratic programming. Technical Report DFVLR-FB 88-28, DLR German Aerospace Center, Institute for Flight Mechanics, Koln, Germany.

\bibitem[{Lee et~al.(2007)Lee, Leok, \protect\BIBand{} McClamroch}]{TaeLeoClam2007spacecraft}
Lee T, Leok M, McClamroch NH (2007) A combinatorial optimal control problem for spacecraft formation reconfiguration. \emph{2007 46th IEEE Conference on Decision and Control} (IEEE), \urlprefix\url{http://dx.doi.org/10.1109/cdc.2007.4434143}.

\bibitem[{L{\'o}pez-Ib{\'a}{\~n}ez et~al.(2022)L{\'o}pez-Ib{\'a}{\~n}ez, Chicano, \protect\BIBand{} Gil-Merino}]{LopChiGil2022evo}
L{\'o}pez-Ib{\'a}{\~n}ez M, Chicano F, Gil-Merino R (2022) The asteroid routing problem: A benchmark for expensive black-box permutation optimization. Jim{\'e}nez~Laredo JL, et~al., eds., \emph{EvoApplications 2022: Applications of Evolutionary Computation}, volume 13224 of \emph{Lecture Notes in Computer Science}, 124--140 (Switzerland: Springer Nature), \urlprefix\url{http://dx.doi.org/10.1007/978-3-031-02462-7_9}.

\bibitem[{Nocedal \protect\BIBand{} Wright(2006)}]{NocWri2006}
Nocedal J, Wright SJ (2006) \emph{Numerical Optimization}. Springer Series in Operations Research and Financial Engineering (Springer), 2nd edition.

\bibitem[{Ow \protect\BIBand{} Morton(1988)}]{OwMor1988:ijpr}
Ow PS, Morton TE (1988) Filtered beam search in scheduling. \emph{International Journal of Production Research} 26:297--307.

\bibitem[{Perez \protect\BIBand{} Régin(2018)}]{Perez_Regin_2018}
Perez G, Régin JC (2018) Parallel algorithms for operations on multi-valued decision diagrams. \emph{Proceedings of the AAAI Conference on Artificial Intelligence} 32(1).

\bibitem[{Petropoulos et~al.(2014)Petropoulos, Bonfiglio, Grebow, Lam, Parker, Arrieta, Landau, Anderson, Gustafson, Whiffen, Finlayson, \protect\BIBand{} Sims}]{PetBonGre2014gtoc5}
Petropoulos AE, Bonfiglio EP, Grebow DJ, Lam T, Parker JS, Arrieta J, Landau DF, Anderson RL, Gustafson ED, Whiffen GJ, Finlayson PA, Sims JA (2014) {GTOC5}: Results from jet propulsion laboratory. \emph{Acta Futura} 8:21--27, \urlprefix\url{http://dx.doi.org/10.2420/AF08.2014.21}.

\bibitem[{Rudich et~al.(2022)Rudich, Cappart, \protect\BIBand{} Rousseau}]{RudCapRou2022cp}
Rudich I, Cappart Q, Rousseau LM (2022) {Peel-And-Bound}: Generating stronger relaxed bounds with multivalued decision diagrams. Solnon C, ed., \emph{Principles and Practice of Constraint Programming}, volume 235 of \emph{LIPIcs}, 35:1--35:20 (Schloss Dagstuhl -- Leibniz-Zentrum f{\"u}r Informatik, Germany), ISBN 978-3-95977-240-2, \urlprefix\url{http://dx.doi.org/10.4230/LIPIcs.CP.2022.35}.

\bibitem[{Rudich et~al.(2023)Rudich, Cappart, \protect\BIBand{} Rousseau}]{RudCapRou2023improved_pnb}
Rudich I, Cappart Q, Rousseau LM (2023) Improved {Peel}-and-{Bound}: Methods for generating dual bounds with multivalued decision diagrams. \emph{Journal of Artificial Intelligence Research} 77:1489--1538, \urlprefix\url{http://dx.doi.org/10.1613/jair.1.14607}.

\bibitem[{Santucci \protect\BIBand{} Baioletti(2022)}]{SanBai2022permutations}
Santucci V, Baioletti M (2022) A fast randomized local search for low budget optimization in black-box permutation problems. \emph{Proceedings of the 2022 World Congress on Computational Intelligence (WCCI 2022)} (Piscataway, NJ: IEEE Press).

\bibitem[{Shirazi et~al.(2018)Shirazi, Ceberio, \protect\BIBand{} Lozano}]{ShiCebLoz2018space}
Shirazi A, Ceberio J, Lozano JA (2018) Spacecraft trajectory optimization: A review of models, objectives, approaches and solutions. \emph{Progress in Aerospace Sciences} 102:76--98, \urlprefix\url{http://dx.doi.org/10.1016/j.paerosci.2018.07.007}.

\bibitem[{Sim{\~o}es et~al.(2017)Sim{\~o}es, Izzo, Haasdijk, \protect\BIBand{} Eiben}]{SimIzzHaas2017multi}
Sim{\~o}es LF, Izzo D, Haasdijk E, Eiben AE (2017) Multi-rendezvous spacecraft trajectory optimization with {Beam} {P-ACO}. Hu B, L{\'o}pez-Ib{\'a}{\~n}ez M, eds., \emph{Proceedings of EvoCOP 2017 -- 17th European Conference on Evolutionary Computation in Combinatorial Optimization}, volume 10197 of \emph{Lecture Notes in Computer Science}, 141--156 (Heidelberg, Germany: Springer), \urlprefix\url{http://dx.doi.org/10.1007/978-3-319-55453-2}.

\bibitem[{Vasile \protect\BIBand{} De~Pascale(2006)}]{VasPas2006preliminary}
Vasile M, De~Pascale P (2006) Preliminary design of multiple gravity-assist trajectories. \emph{Journal of Spacecraft and Rockets} 43(4):794--805, \urlprefix\url{http://dx.doi.org/10.2514/1.17413}.

\bibitem[{Zaefferer et~al.(2014)Zaefferer, Stork, Friese, Fischbach, Naujoks, \protect\BIBand{} Bartz-Beielstein}]{ZaeStoFriFisNauBar2014}
Zaefferer M, Stork J, Friese M, Fischbach A, Naujoks B, Bartz-Beielstein T (2014) Efficient global optimization for combinatorial problems. Igel C, Arnold DV, eds., \emph{Proceedings of the Genetic and Evolutionary Computation Conference, GECCO 2014}, 871--878 (New York, NY: ACM Press), \urlprefix\url{http://dx.doi.org/10.1145/2576768.2598282}.

\end{thebibliography}

\end{document}